\documentclass[11pt]{article}

\usepackage{bridges}
\usepackage{amsfonts,amssymb,amsthm,eucal,amsmath}
\usepackage{graphicx}
\usepackage{subfig}

\usepackage{wrapfig}
\usepackage{pinlabel, color}
\usepackage{multirow}

\usepackage{microtype}
\usepackage{hyperref}
\usepackage[para]{manyfoot}
\DeclareNewFootnote[para]{B}

\newcommand{\EE}{\mathbb{E}}
\newcommand{\HH}{\mathbb{H}}

\newcommand{\RR}{\mathbb{R}}

\newcommand{\arcsinh}{\text{arcsinh}}

\newcommand{\HxE}{\mathbb{H}^2\times\mathbb{E}}

\newcommand{\black}{\color{black}}

\begin{document}

\title{Non-euclidean virtual reality II: explorations of $\HxE$}
\author{
 \begin{tabular}{cccc}
  Vi Hart & Andrea Hawksley & Elisabetta A. Matsumoto & Henry Segerman \\
  eleVR & eleVR & School of Physics & Department of Mathematics  \\
  HARC & HARC & Georgia Institute of Technology & Oklahoma State University 
 \end{tabular}
}
\date{}

\maketitle

\begin{abstract}
We describe our initial explorations in simulating non-euclidean geometries in virtual reality. Our simulation of the product of two-dimensional hyperbolic space with one-dimensional euclidean space is available at \href{http://h2xe.hypernom.com}{h2xe.hypernom.com}.\footnote{The code is available at \href{https://github.com/henryseg/H2xE_VR}{github.com/henryseg/H2xE\_VR}.} 
\end{abstract}

\begin{figure}[htb]
\centering
{
\includegraphics[width=0.9\textwidth]{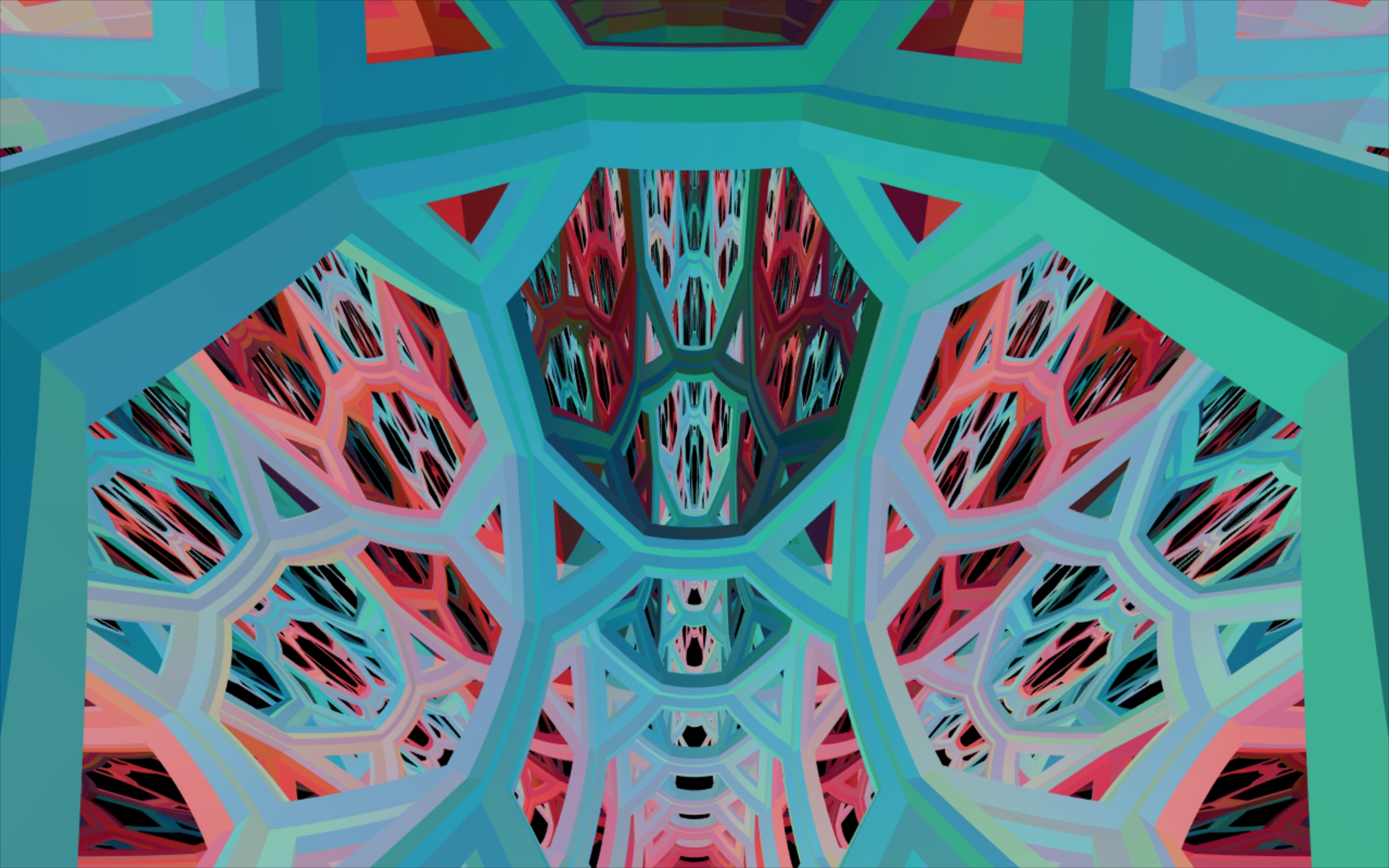}
\label{h2xr_title}
}
\caption{A view from $\HH^2\times \EE$.}
\label{H2xE_title}
\end{figure}

The properties of euclidean space seem natural and obvious to us, to the point that it took mathematicians over two thousand years to see an alternative to Euclid's parallel postulate. The eventual discovery of hyperbolic geometry in the 19th century shook our assumptions, revealing just how strongly our native experience of the world blinded us from consistent alternatives, even in a field that many see as purely theoretical. Non-euclidean spaces are still seen as unintuitive and exotic, but we believe that with direct immersive experience we can get a better ``feel" for them. The latest wave of virtual reality hardware, in particular the HTC Vive, tracks both the orientation and the position of the headset within a room-sized volume, allowing for such an experience.

Most visualizations of hyperbolic space are seen from the outside, as in Escher's \emph{Circle Limit} series of prints, which use the Poincar{\'e} disk model of two-dimensional hyperbolic space. 
Three-dimensional hyperbolic space can also be visualised
in a similar way, via the Poincar\'e ball model. In virtual reality, we could simulate this floating in the middle of the room. 
We would then generate graphics on screen using the standard euclidean graphics pipeline, and motions in real-life would translate directly to motions of the virtual camera in the 
ambient euclidean space of the Poincar\'e ball model.
Even though you would be able to put your head inside of this virtual Poincar\'e ball model, it would give an ``extrinsic'' experience of $\HH^3$ -- you would experience the induced metric of the Poincar\'e model and not directly the metric of $\HH^3$.
Such extrinsic visualisations provide a brief, compact snapshot of infinite hyperbolic space, yet the viewer is left to their own imagination to remodel the space to see what life might be like for an inhabitant living inside. 
Our goal is to make three-dimensional non-euclidean spaces feel more natural by giving people experiences inside those spaces, including the ability to move through those spaces with their bodies. 
Luckily, computers don't know or care that people live in a mostly euclidean world, a world where cubes fit together four around an edge because they have $90^\circ$ angles. As long as we program in the correct mathematics, a computer is perfectly happy simulating a hyperbolic space where cubes pack neatly, six around an edge.

\begin{figure}[t]
\vspace{-20pt}\centering
\subfloat[View in the $\HH^2$ direction.]
{
\includegraphics[width=0.45\textwidth]{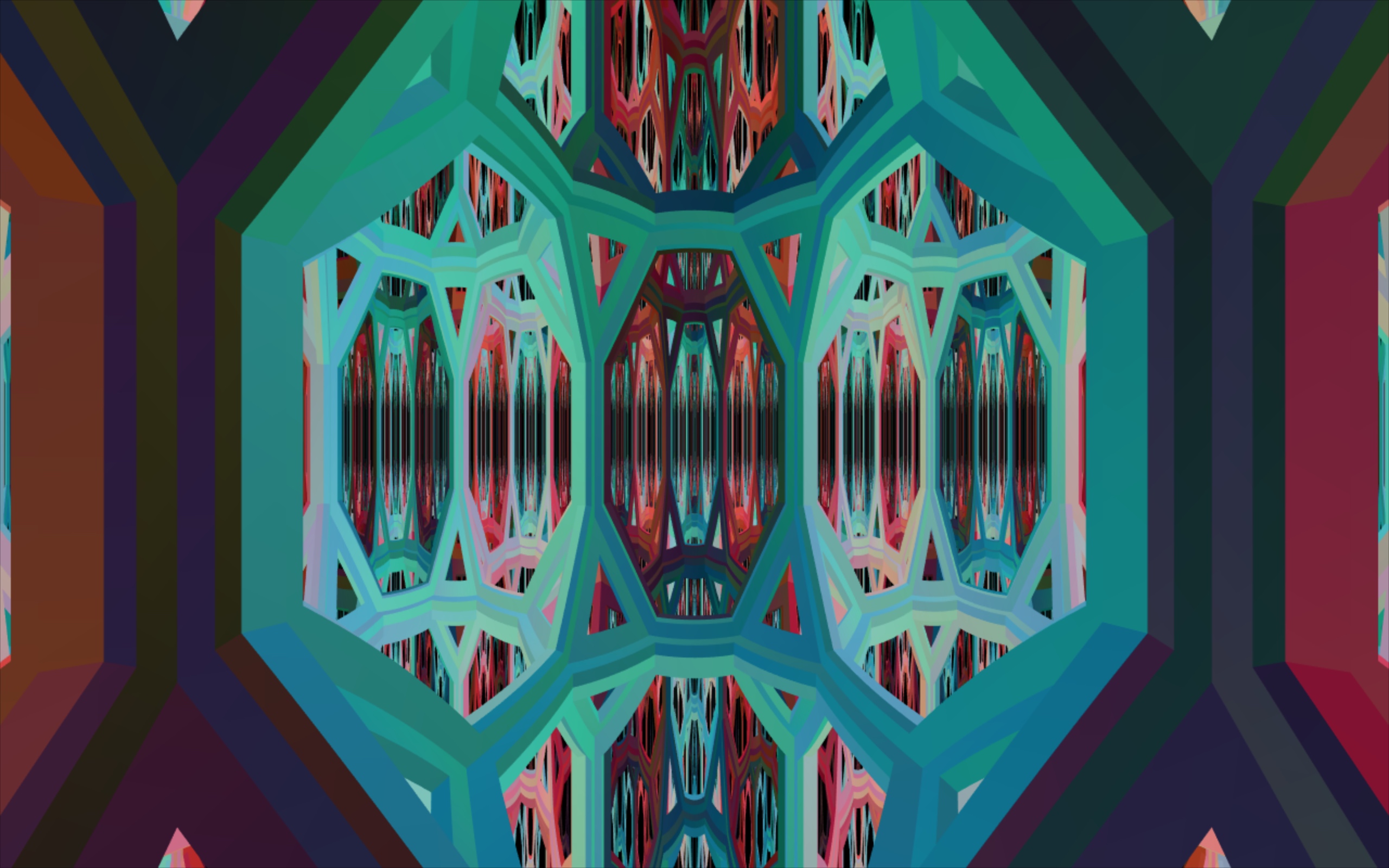}
\label{h2xr_46_a}
}
\quad
\subfloat[View in a diagonal direction.]
{
\includegraphics[width=0.45\textwidth]{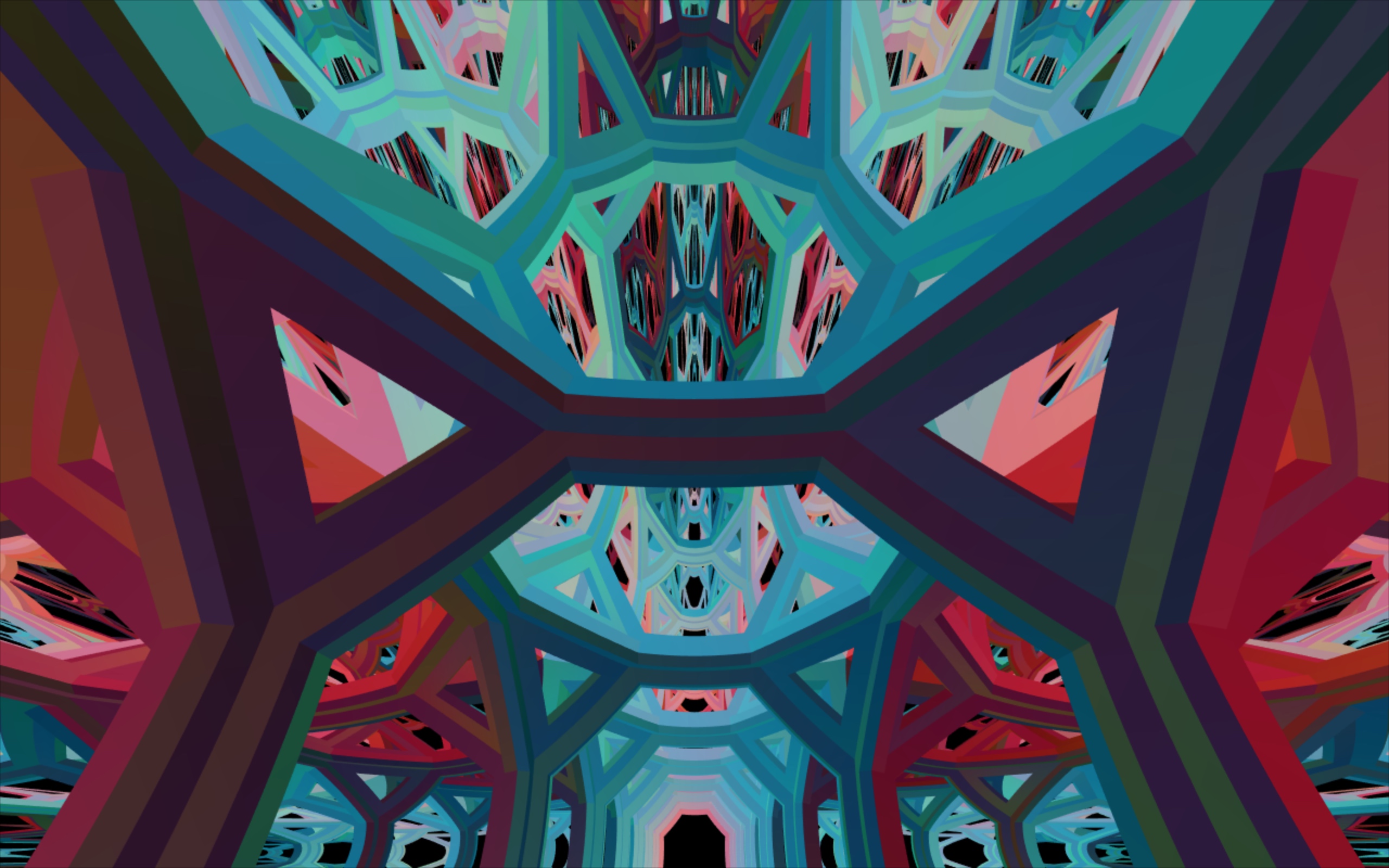}
\label{h2xr_46_b}
}

\vspace{-5pt}
\subfloat[View in the $\EE$ direction.]
{
\includegraphics[width=0.45\textwidth]{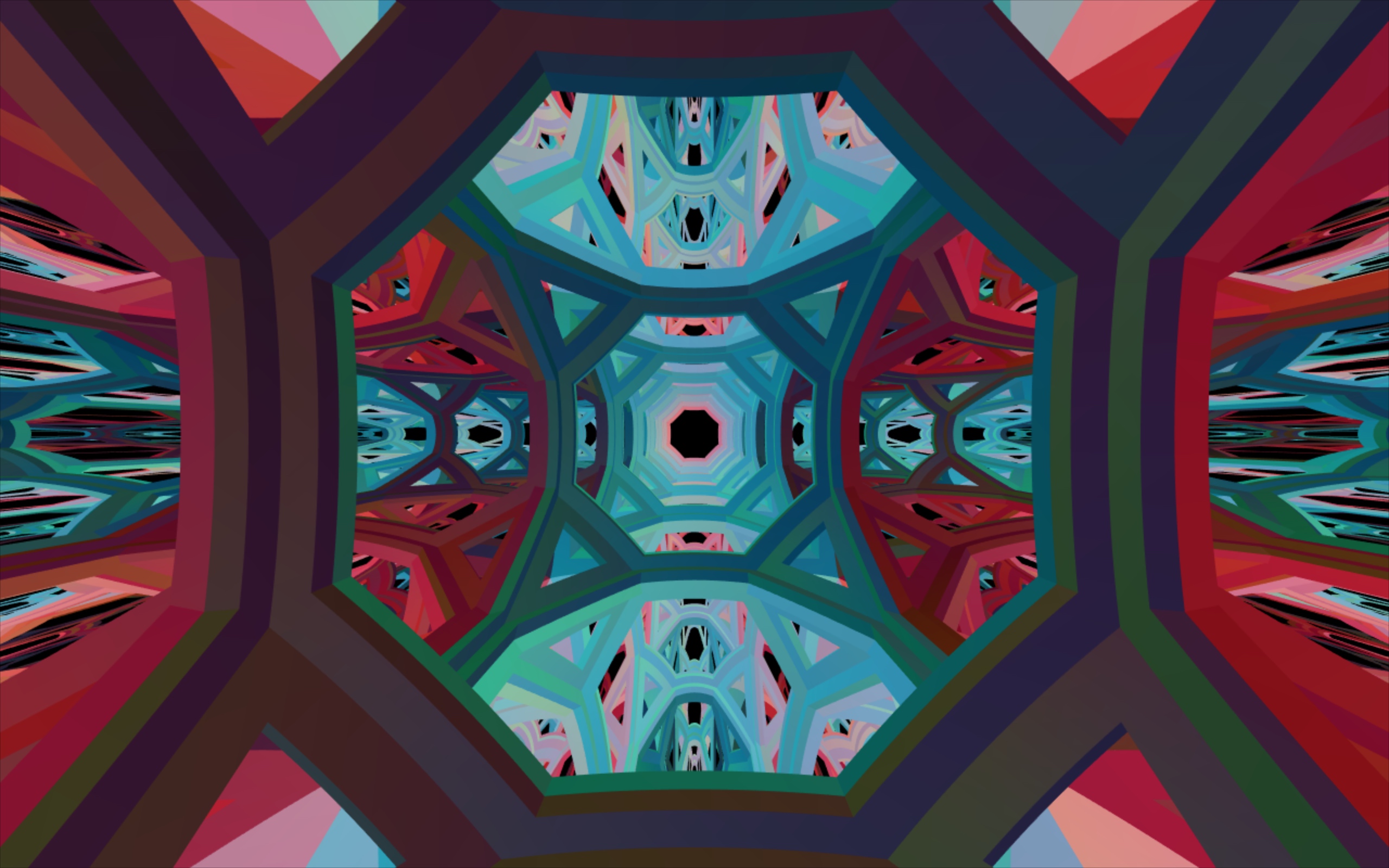}
\label{h2xr_46_c}
}
\quad
\subfloat[Another view in a diagonal direction.]
{
\includegraphics[width=0.45\textwidth]{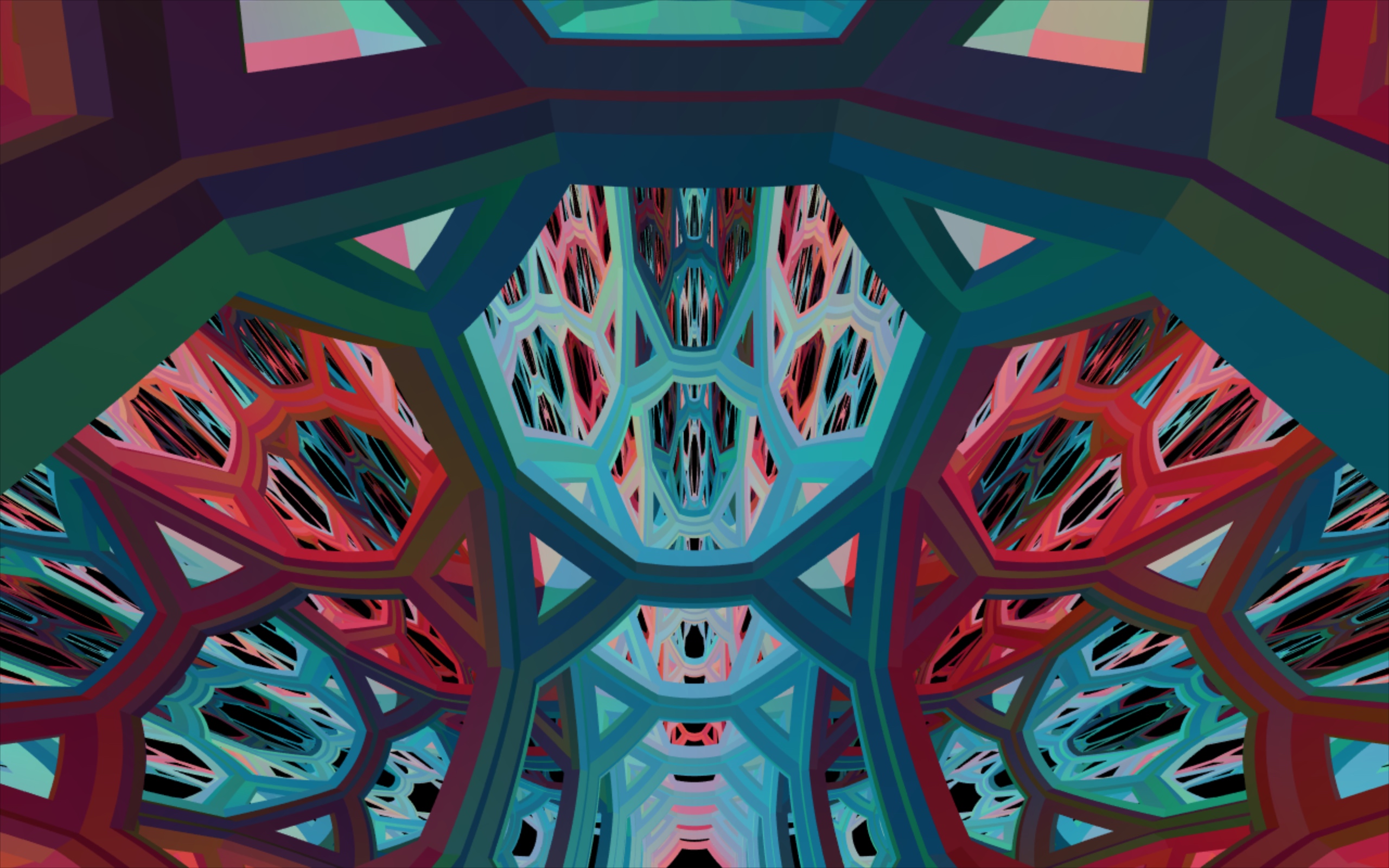}
\label{h2xr_46_d}
}

\caption{Views of the $\{4,6\}$ tiling of $\HH^2\times\EE$. We draw the honeycomb to a depth of seven steps from the central cube.}
\label{h2xr_46}
\vspace{-15pt}
\end{figure}

We took inspiration from Jeff Weeks' \emph{Curved Spaces}~\cite{curved_spaces} software, a ``flight simulator for multiconnected universes'', which allows the user to ``fly'' a spaceship through various three-dimensional manifolds, with $S^3, \EE^3$ and $\HH^3$ geometries. The user controls the direction of the spaceship using the mouse, and its speed using keyboard controls. With virtual reality technology however, the user controls the direction in which they are looking by turning their head, and their position by moving their body. 
Thus, we remove barriers between us and the experience of the space -- it is easier and more natural to access the experience, particularly for users who are not familiar with moving through space using ``computer game'' controls, and this extra ease allows a user to discover some not-so-obvious properties of these spaces much more readily.

Spaces such as $\HH^3$ and $S^3$ are appealing to create and explore because they are \emph{isotropic}, a property shared with euclidean space -- the space is uniform in any direction you look or walk. Yet, when you move through these spaces, you encounter several unexpected consequences stemming from \emph{parallel transport}, rotations of the reference frame as you traverse a curved space~\cite{our_paper_h3}. One of the more disconcerting experiences in an $\HH^3$ virtual reality simulator is watching the floor fall out from beneath your feet as you walk forward. We sought a gentler introduction into a hyperbolic space. One of the other eight Thurston geometries~\cite{thurston_book}, $\HxE$ -- the cartesian product of the hyperbolic plane and the euclidean line -- 
enables the user to traverse hyperbolic space horizontally as they walk through a room, yet it retains familiar euclidean geometry in the vertical direction. In this space, the virtual floor remains
in the same place as the real-life floor. Jeff Weeks has also explored $\HH^2 \times \EE$ via a similar method to ours in unpublished work.

There are four ingredients that go into our virtual reality simulation of $\HxE$ as outlined in this paper:
\begin{enumerate}
\item{A way to describe the points of $\HxE$ numerically, i.e. a \emph{model} of $\HxE$}

\item{A way to convert points in the model into points in $\EE^3$ that we can then draw on screen,}

\item{A way to move around $\HxE$ using the motion inputs from the virtual reality headset, and}

\item{A set of landmarks in $\HxE$ to draw, to help the viewer navigate the space -- 
we use a tiling of $\HxE$.}

\end{enumerate}

\section{The Model of $\HxE$}
The space $\HH^2 \times \EE$ is the cartesian product of the two topological spaces $\HH^2$ and $\EE$. Unlike $\EE^3$, $S^3$ and $\HH^3$, $\HH^2 \times \EE$ is not isotropic: the geometry is different when we look along each of the $\HH^2$ and the $\EE$ directions. 
Our model is the cartesian product of the hyperboloid model of $\HH^2$ with the real line. The model of $\HH^2\times\EE$ lives in four-dimensional Minkowski space $\EE^{3,1}$ as the set of points $\{(x,y,z,w) \in \EE^{3,1} \mid x^2+y^2 = w^2 - 1, w>0\}.$  The coordinates $x,y$ and $w$ are used for the hyperboloid model of $\HH^2$, while the coordinate $z$ is used for $\EE$. We use the explicit parametrisation given by the map 
$\boldsymbol{\phi}: \RR_{\ge 0} \times  [0,2\pi) \times \RR \rightarrow \EE^{3,1}$ given by 
\[
\boldsymbol{\phi}(\rho,\theta,z) = (\sinh\rho\cos\theta,\sinh\rho\sin\theta,z,\cosh\rho).
\]

\section{Drawing points in $\HxE$ on the screen}

\begin{wrapfigure}[12]{r}{0.27\textwidth}
  \vspace{-15pt}
  \centering
  \includegraphics[width=0.22\textwidth]{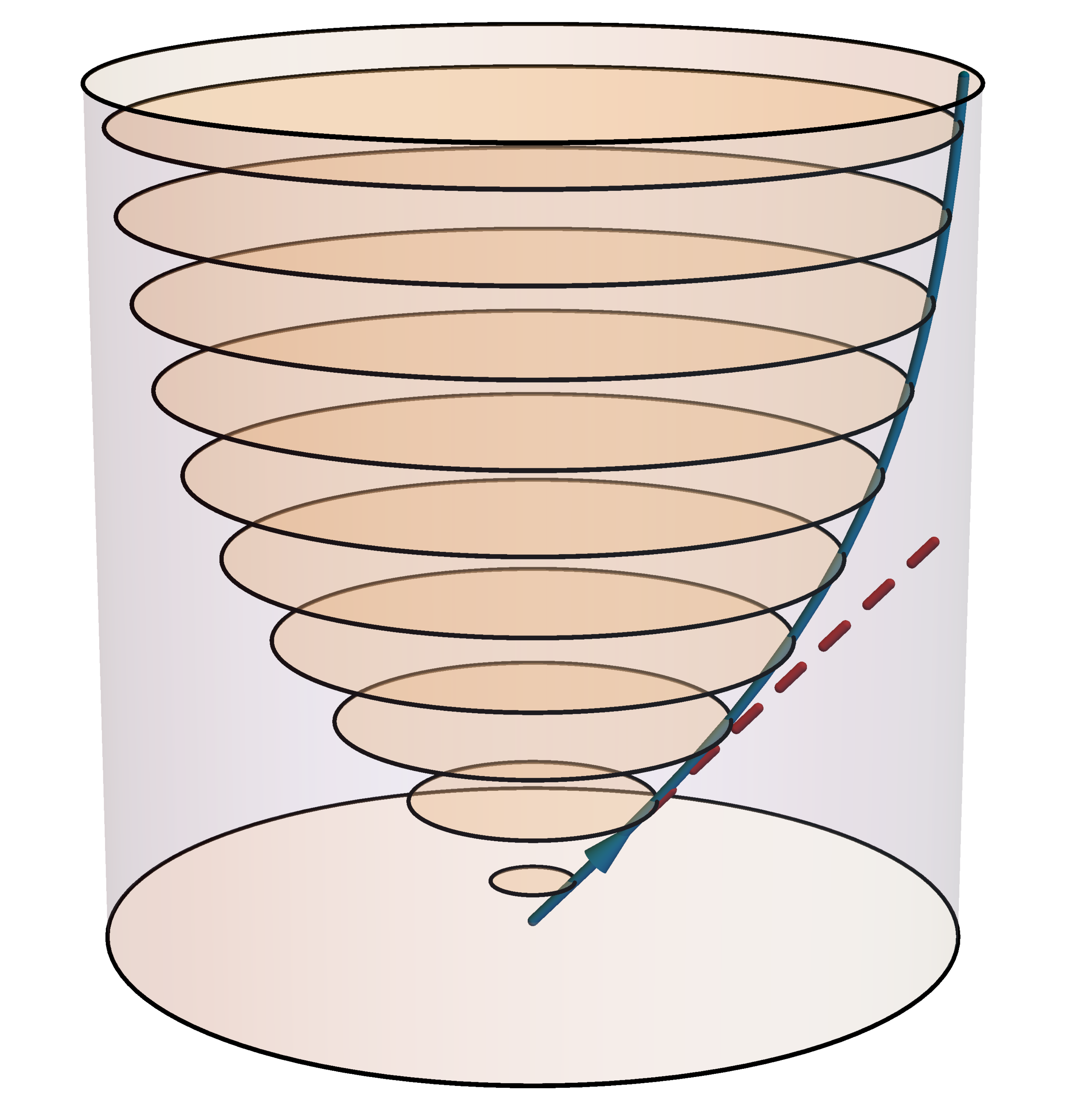}
  \caption{A geodesic in the product of the Klein model with $\EE$, and a straight line in $\EE^3$.} 
  \label{Fig:Klein_x_R_geodesic}
\end{wrapfigure}
To draw points on screen, we cannot use the same trick as we did in $\HH^3$ of ignoring the inverse of the exponential map and using the Klein model instead~\cite{our_paper_h3}. Na\"ively, we might try converting the $x,y$ and $w$ coordinates of a point into a point in the Klein model of $\HH^2$, sitting inside of the unit disk in $\RR^2$, and draw points in the cartesian product of the Klein model and $\EE$, using $z$ for the value in $\EE$. If we did this then we would see objects in the space, but if we tried to move towards an object we would miss. This is because the light rays in $\EE^3$ which connect an object to our eyes do not lie along the geodesics of (Klein)$\times\EE$. To properly account for this, we shall digress to a discussion of geodesics and the exponential map in $\HH^2\times \EE.$

Whilst geodesics in the Klein model are straight lines in the euclidean space in which the model lives, geodesics in  (Klein model)$\times \EE$ have an altogether different flavour due to the anisotropic nature of $\HH^2\times\EE$. 
If we are at the origin, looking in the diagonal direction shown in Figure \ref{Fig:Klein_x_R_geodesic}, we see objects that are on the red line. But when we move in that direction we arrive at objects that are on the blue geodesic.

So, we need to return to the original plan in our previous paper~\cite{our_paper_h3} of calculating the inverse of the exponential map in order to draw objects correctly. To 
do this, we need to 
work out parametrisations of geodesics in our model of $\HH^2\times\EE$ in $\EE^{3,1}$.

The first step in calculating geodesics is finding the components of the metric tensor. In our parametrisation $\boldsymbol{\phi}(\rho,\theta,z) = (\sinh\rho\cos\theta,\sinh\rho\sin\theta,z,\cosh\rho)$, the components of the metric tensor $g_{ij}$ are given by $g_{ij}=\langle\partial_i \boldsymbol{\phi},\partial_j\boldsymbol{\phi}\rangle,$ where $i,j \in \{\rho,\theta,z\}$ and $\langle \cdot,\cdot\rangle$ denotes the inner product in the ambient space, $\EE^{3,1}$. The three non-zero components of the metric tensor are $g_{\rho\rho}=g_{zz}=1$ and $g_{\theta\theta}=(\sinh\rho)^{-2}$. Each of the derivatives $\partial_i \boldsymbol{\phi}$ is a vector in the tangent space of $\HH^2\times\EE$.

The derivative of a vector in the tangent space tells us the rate of change of that vector (both magnitude and direction) as we move along another vector in the tangent space.
If the derivative has no component in any direction other than parallel to itself, then parallel transport in the direction of the vector preserves it. This is our condition for being a geodesic. Thus the geodesic $\boldsymbol{\gamma}(t)=\boldsymbol{\phi}(\rho(t),\theta(t),z(t))$ satisfies the equation $\nabla_{\dot{\boldsymbol{\gamma}}}\dot{\boldsymbol{\gamma}}=0$, where $\dot{\boldsymbol{\gamma}}=\frac{d\boldsymbol{\gamma}}{dt}$ and $\nabla_X$ is the Levi-Civita connection. In coordinate form this may be re-written as 
\[
\partial^2_t\boldsymbol{\gamma}^\lambda+\sum_{\mu,\nu}\Gamma^{\lambda}_{\mu\nu}\partial_t\boldsymbol{\gamma}^\mu\partial_t\boldsymbol{\gamma}^\nu=0,
\]
where Christoffel symbols of the second kind $\Gamma^{\lambda}_{\mu\nu} = \frac{1}{2}\sum_{\sigma}g^{-1}_{\lambda\sigma}(\partial_\nu g_{\sigma\mu}+\partial_\mu g_{\sigma\nu}-\partial_\sigma g_{\mu\nu})$ are the components of the Levi-Civita connection and $\lambda,\mu,\nu,\sigma \in \{\rho,\theta,z\}$. In $\HH^2\times\EE$, the only non-zero Christoffel symbols are $\Gamma^{\rho}_{\theta\theta}=-\cosh \rho\sinh\rho$, and $\Gamma^{\theta}_{\rho\theta}=\Gamma^\theta_{\theta\rho}=\coth \rho$. This set of coupled differential equations gives a formula for parametrised curves in the coordinate map $\boldsymbol\gamma(t)=(\rho(t),\theta(t),z(t))$.
In $\HxE$ geodesics emanating from the origin $(0,0,0,1)$ of the hyperboloid must satisfy the differential equations for curves $\rho''(t)=z''(t)=\theta'(t)=0,$ with boundary conditions $\rho(0)=0,$ $z(0)=0$, $\theta(0)=\theta_0$, $\rho(1)=\rho_0$ and $z(1)=z_0$. Thus they are parametric curves in $\EE^{3,1}$ of the form 
\[
\boldsymbol{\gamma}(t) = \{\sinh(\rho_0 t) \cos\theta_0,\sinh(\rho_0 t)\sin\theta_0, z_0 t, \cosh(\rho_0 t)\}.
\]

To draw points on screen, we construct the inverse of the exponential map. Given a point $p$ in $\HH^2\times\EE$, this tells us what initial velocity in the tangent space we need in order to reach $p$ in $\HH^2\times\EE$ at time $t=1$. 
 The inverse of the exponential map sends points on the (hyperboloid)$\times\EE$ in $\EE^{3,1}$ to points in $\EE^3$ via
\[
(x,y,z,w)\rightarrow\left(\frac{\arcsinh(\sqrt{x^2+y^2}) }{\sqrt{x^2+y^2}}x,\frac{\arcsinh(\sqrt{x^2+y^2}) }{\sqrt{x^2+y^2}}y,z\right).
\]
The images of geodesics in $\HH^2 \times \EE$ are straight lines in the tangent space. Thus if you look at a point in the virtual reality experience and move towards it, you will eventually reach it.

\section{Moving through a curved space in virtual reality}

Each frame of our simulation, the virtual reality headset (the HTC Vive), outputs both position and orientation information corresponding to the location of your head in the room and the direction in which you are looking. We can use this information in several ways, depending on what kind of experience we wish to create.\footnote{In the vast majority of virtual reality experiences, the position and orientation of the headset are mapped directly to the position and orientation of a virtual camera in euclidean space. In our geometry simulations, we map orientation directly but treat position differently. In our previous work \emph{Hypernom}~\cite{hypernom}, there is no position tracking, but we map the headset's orientation to both the orientation and the position of the virtual camera in $S^3$.}

There is freedom in how to map the headset data of the user's motion in the room into the hyperboloid model. One possibility would be to map the position and orientation data onto the Poincar\'e ball model, and 
show the in-space view of $\HH^3$ from that position and orientation. Although, this might not be the most obvious choice, it would allow us to 
 map an infinite space into the finite confines of a room in $\EE^3$. Unfortunately, this sort of mapping violates a property of movement that humans are quite attached to, which is consistency in distance. The user would find that two motions of their head of the same magnitude would translate to vastly different translations within the simulation, depending on where in the room they are standing.

\black
The approach we take here is to look at the relative motion of the headset, and move the virtual camera in a corresponding way. Every frame, we compare the headset's current position to its previous position, compute the difference, and move the virtual position by that vector. This has the advantage that your head's motions behave the same way no matter where you are. The disadvantage is that at any one time, any location in your real space might map to any location in the virtual space. 

Once we know the relationship between the viewer's motions and the headset output, we must convert them to motions in $\HxE$. Moving through this space now consists of translating by isometries inherited from $\EE^{3,1}$ that leave the model unchanged.
Turning this set of isometries into the exponential map can be handled in a variety of ways. As both $\HH^2$ and $\EE$ can each be embedded in a flat space of one higher dimension, the exponential map could be encoded in a block diagonal $5\times5$ matrix, yet graphics cards are not optimised for this type of matrix multiplication. Thus, we can take advantage of $\HxE$ being a product space where the isometries of $\EE$ (adding real numbers together gives pure translations along $\EE$) and the isometries of $\HH^2$ (implemented in exactly the same way as we did for $\HH^3$~\cite{our_paper_h3}) act independently of each other.

\section{The $\{4,6\}$ tiling of $\HH^2$ and a colouring of $\HH^2\times\EE$}

\begin{figure}[t]
\centering
\vspace{-10pt}
\subfloat[Our colouring on the $\{4,6\}$ tiling in the Poincar\'e disk model.]
{
\includegraphics[width=0.3\textwidth]{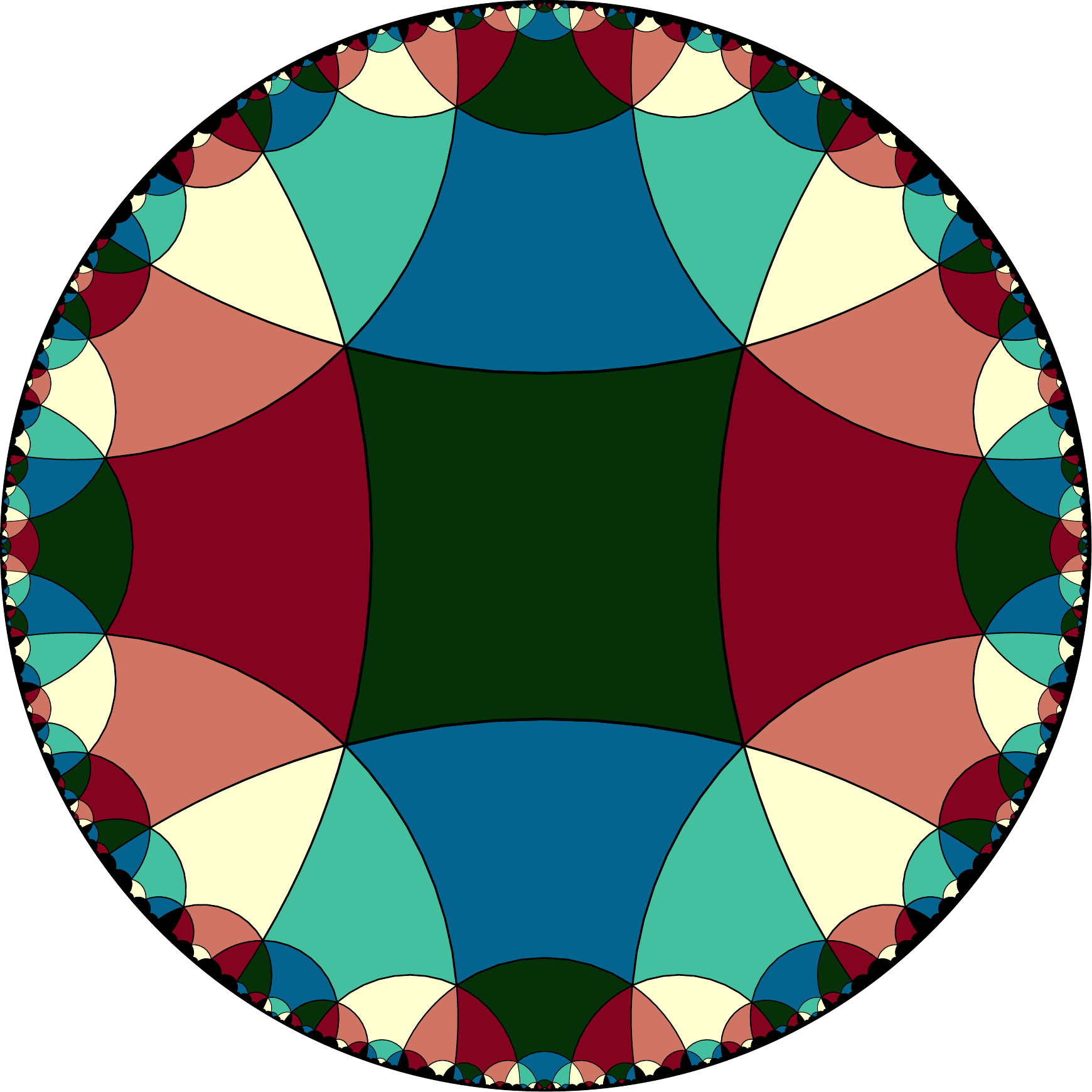}
\label{46_colouring}
}
\quad
\subfloat[A fundamental domain for the tiling, drawn in the Poincar\'e disk model.]
{
\includegraphics[width=0.3\textwidth]{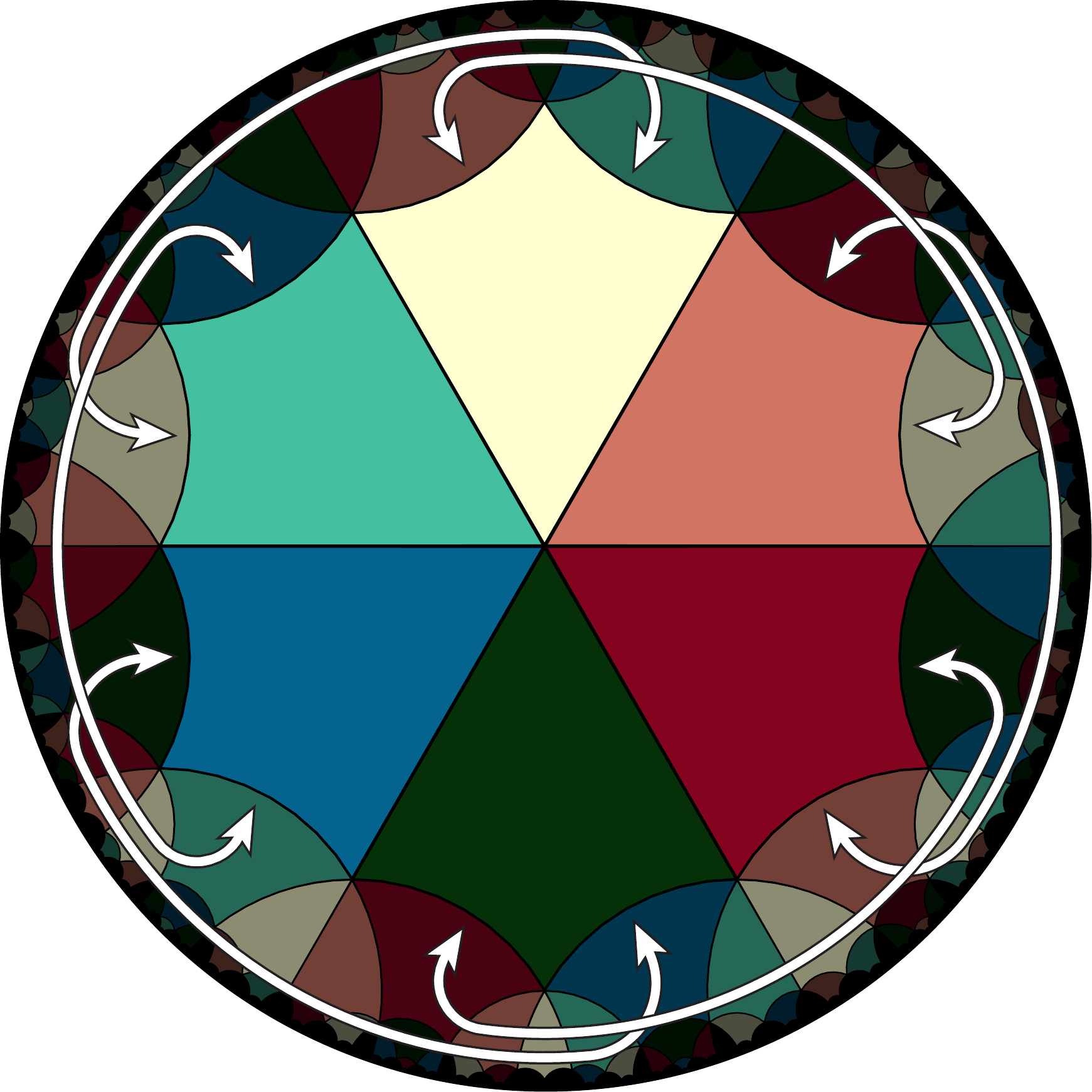}
\label{46_colouring_genus_2_universal_cover}
}
\quad
\subfloat[The fundamental domain glued up to form a genus two surface.]
{
\includegraphics[width=0.3\textwidth]{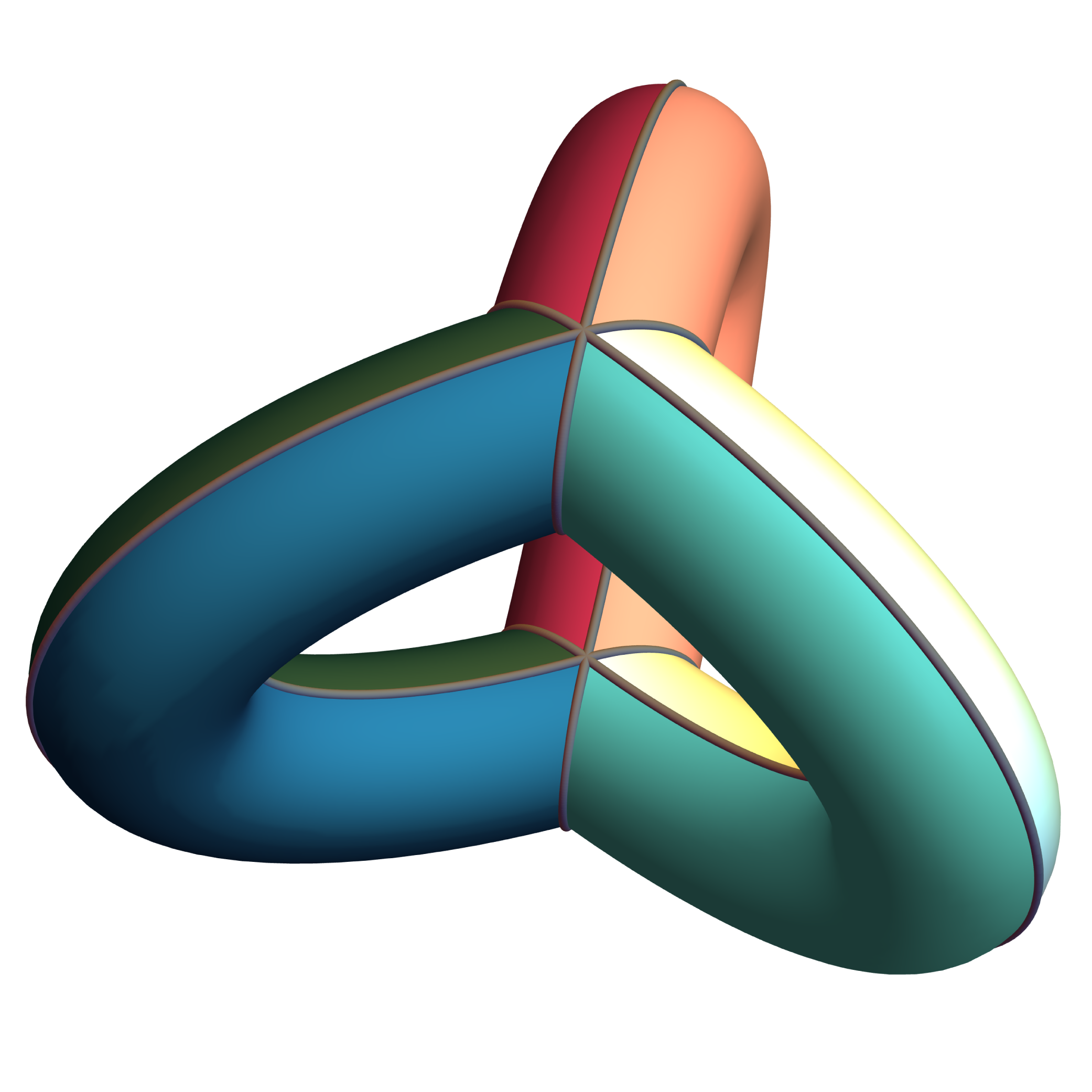}
\label{46_colouring_genus_2}
}
\caption{Our colouring of the $\{4,6\}$ tiling.}
\label{46_colouring}
\vspace{-10pt}
\end{figure}

\begin{figure}[b]
\vspace{-10pt}
\centering
\subfloat[$\HxE$ initial view.]
{
\includegraphics[width=0.3\textwidth]{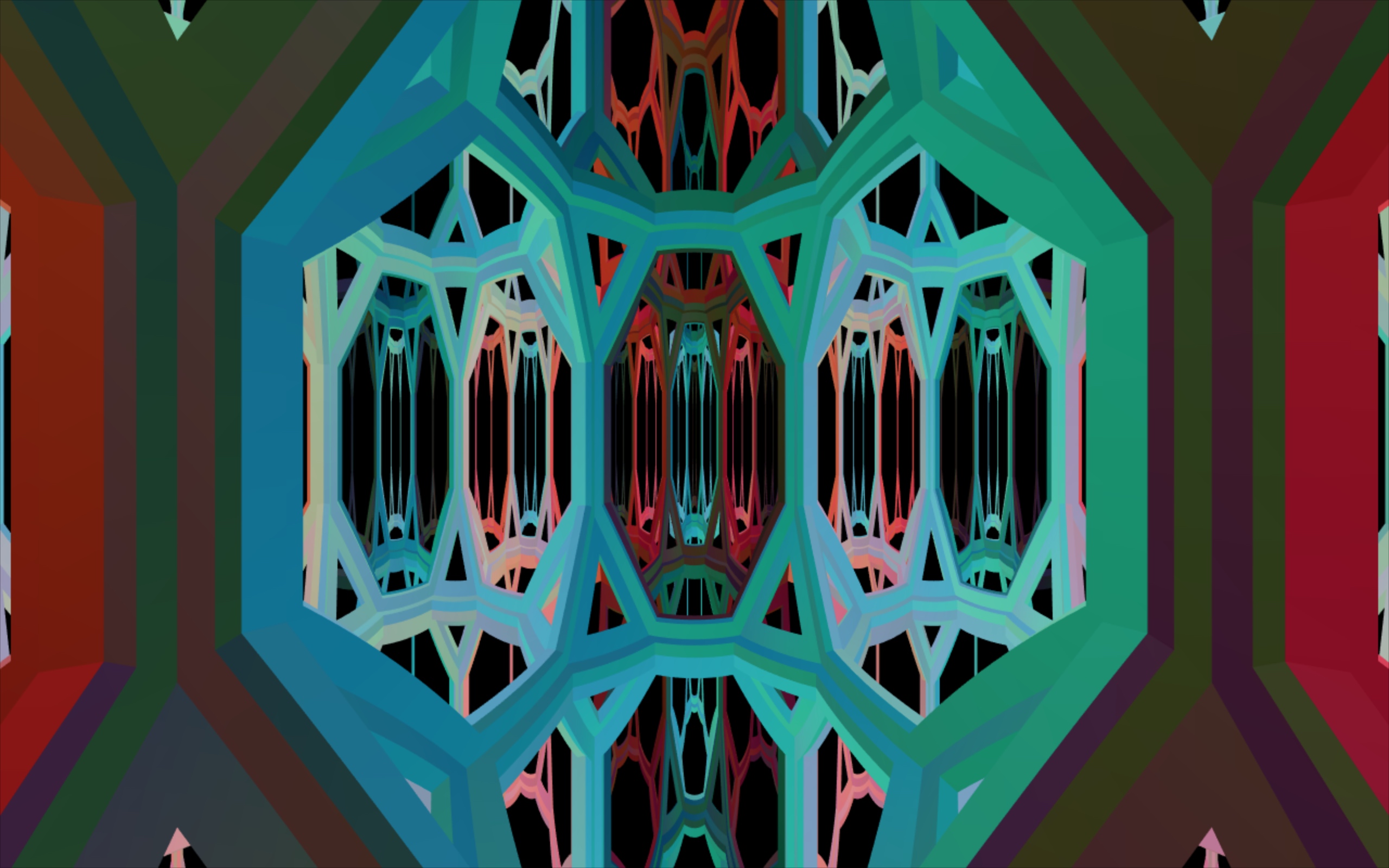}
\label{aspect_ratio_1}
}
\quad
\subfloat[Moving forward 1/6 a cube width.]
{
\includegraphics[width=0.3\textwidth]{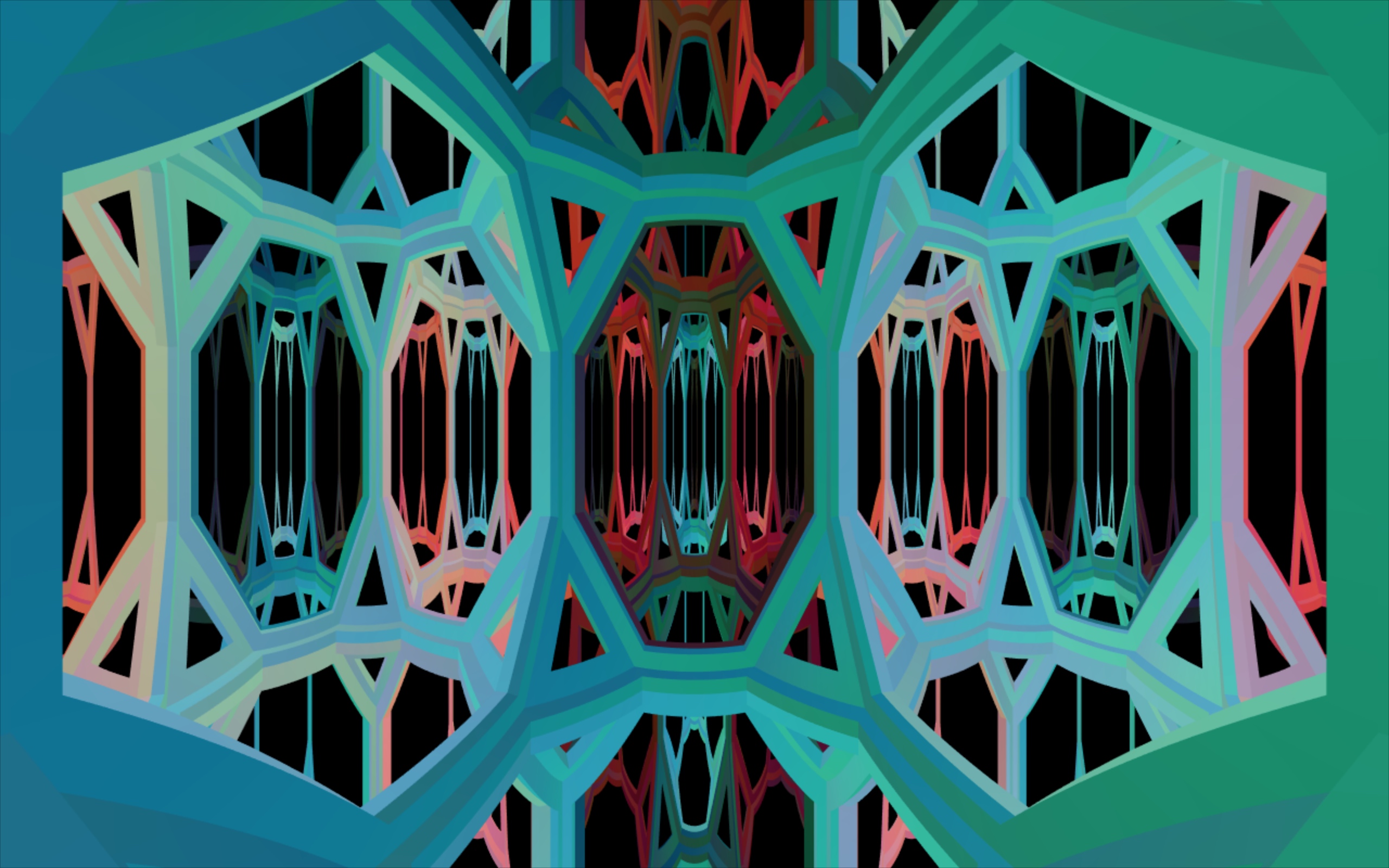}
\label{aspect_ratio_2}
}
\quad
\subfloat[Moving forward 1/3 a cube width.]
{
\includegraphics[width=0.3\textwidth]{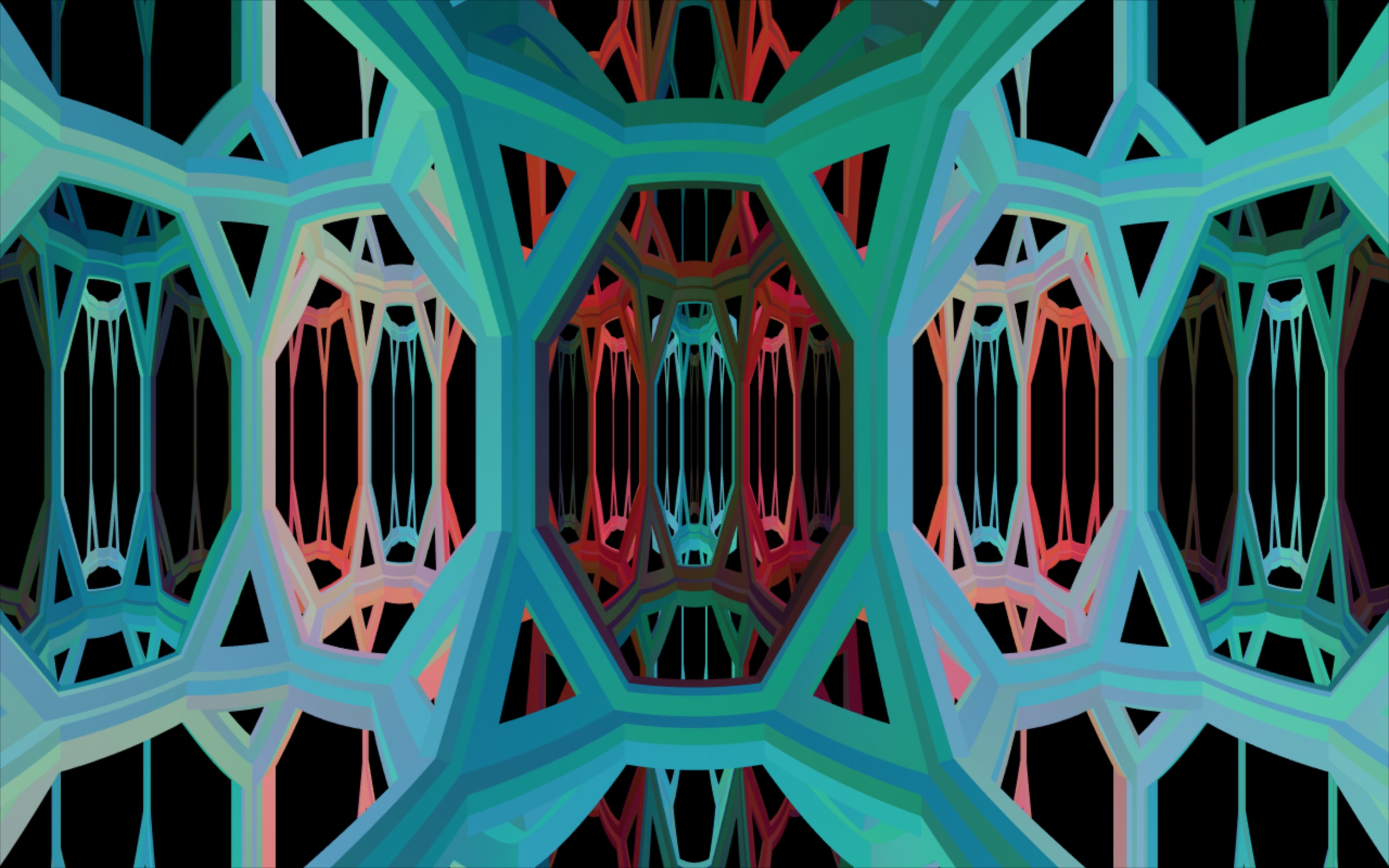}
\label{aspect_ratio_3}
}

\vspace{-7pt}
\subfloat[Moving forward 1/2 a cube width.]
{
\includegraphics[width=0.3\textwidth]{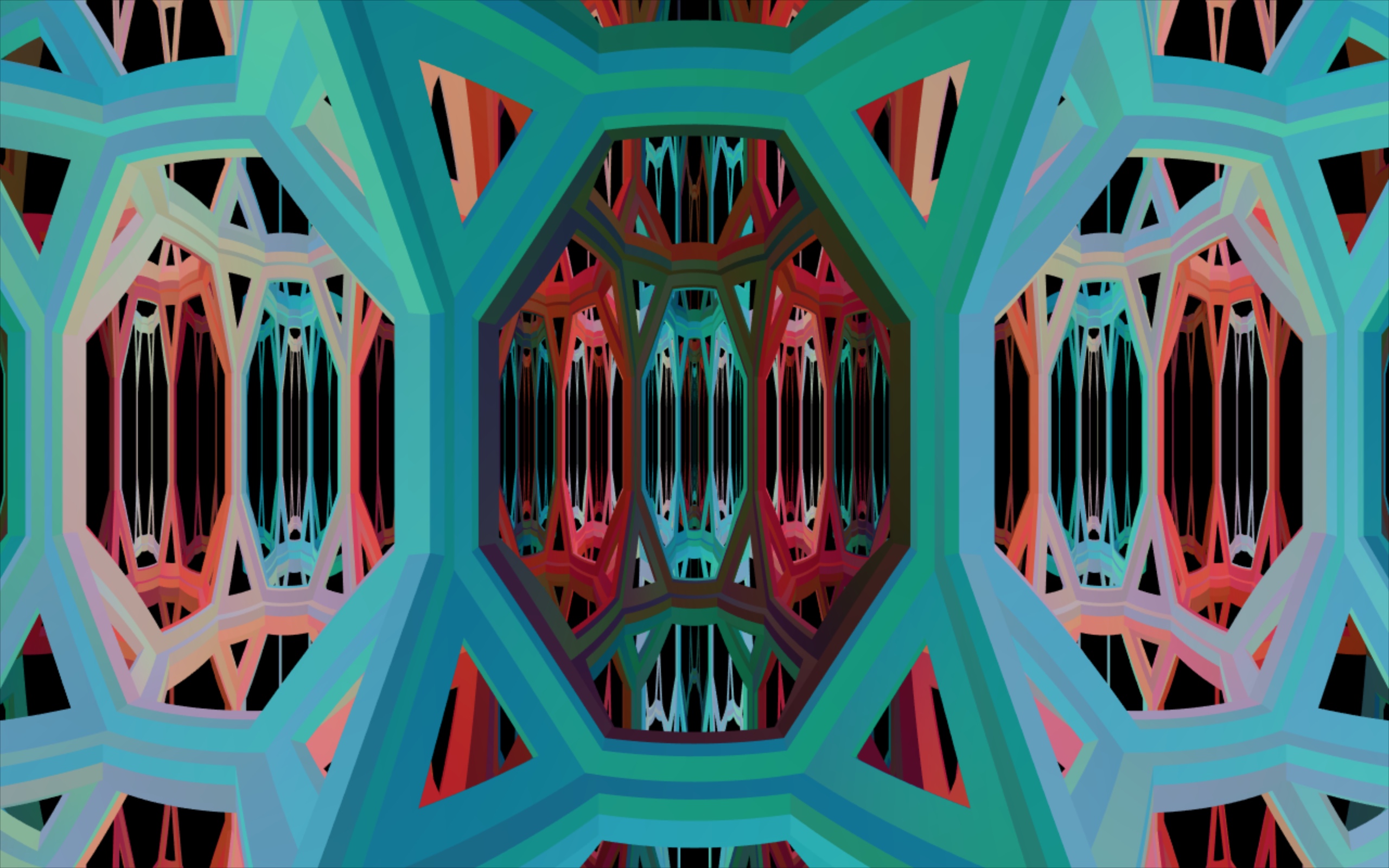}
\label{aspect_ratio_4}
}
\quad
\subfloat[Moving forward 2/3 cube width.]
{
\includegraphics[width=0.3\textwidth]{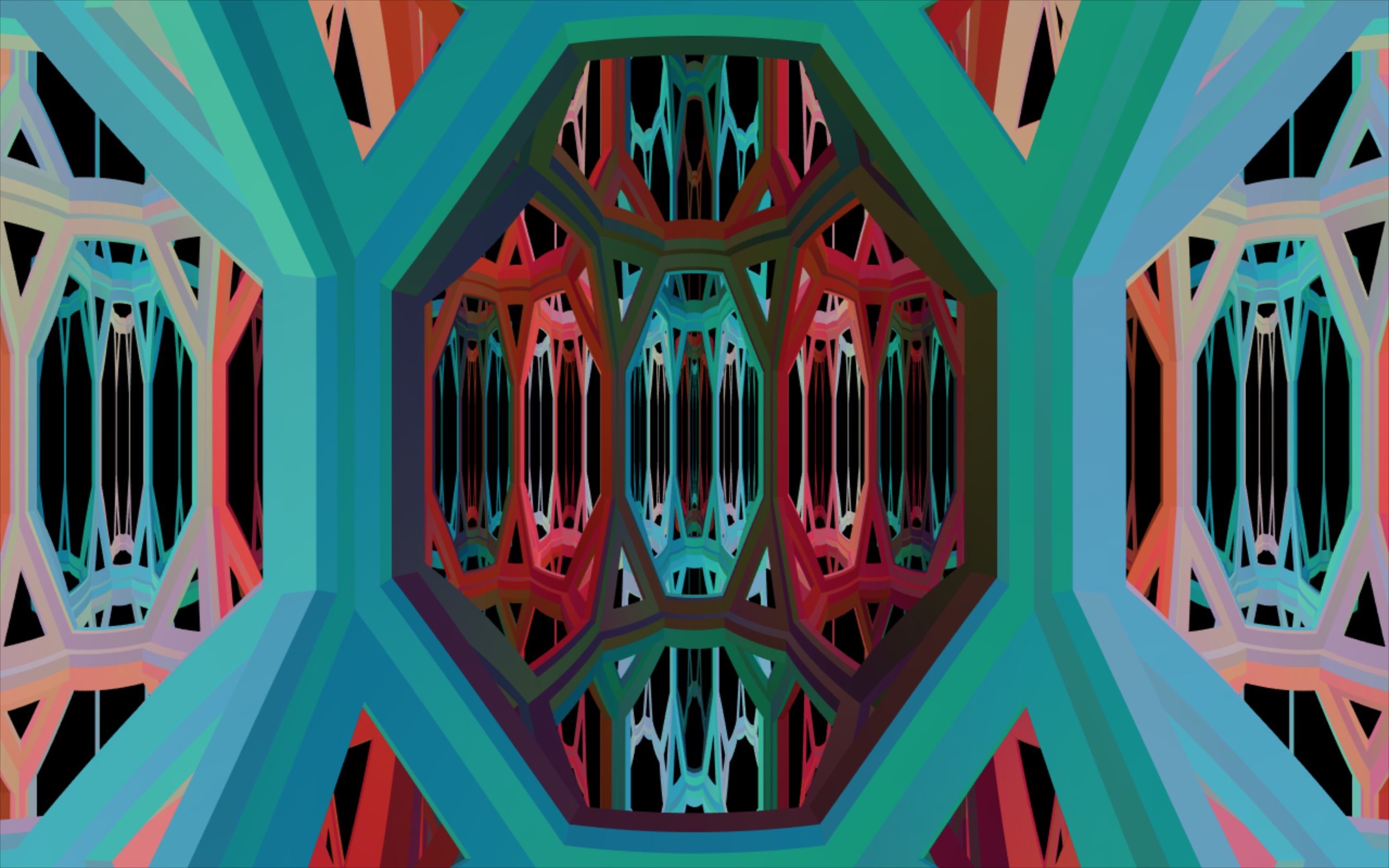}
\label{aspect_ratio_5}
}
\quad
\subfloat[Moving forward 1 cube width.]
{
\includegraphics[width=0.3\textwidth]{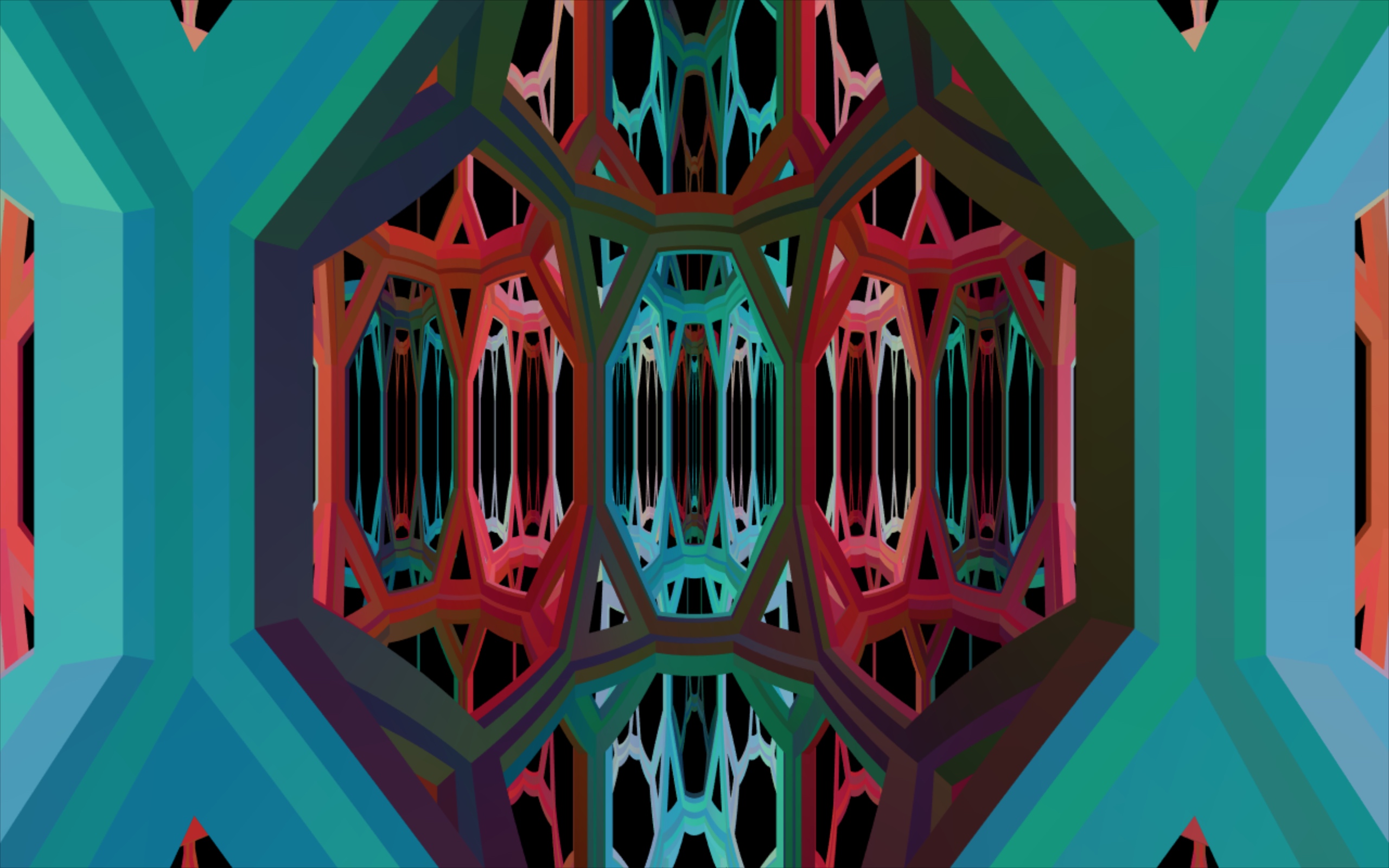}
\label{aspect_ratio_6}
}
\caption{The aspect ratio of the cubes scales differently in the $\HH^2$ and $\EE$ directions.}
\label{Fig:aspect_ratio}
\end{figure}

\begin{figure}[t]
\centering
\subfloat[$\HH^3$ initial view.]
{
\includegraphics[width=0.22\textwidth]{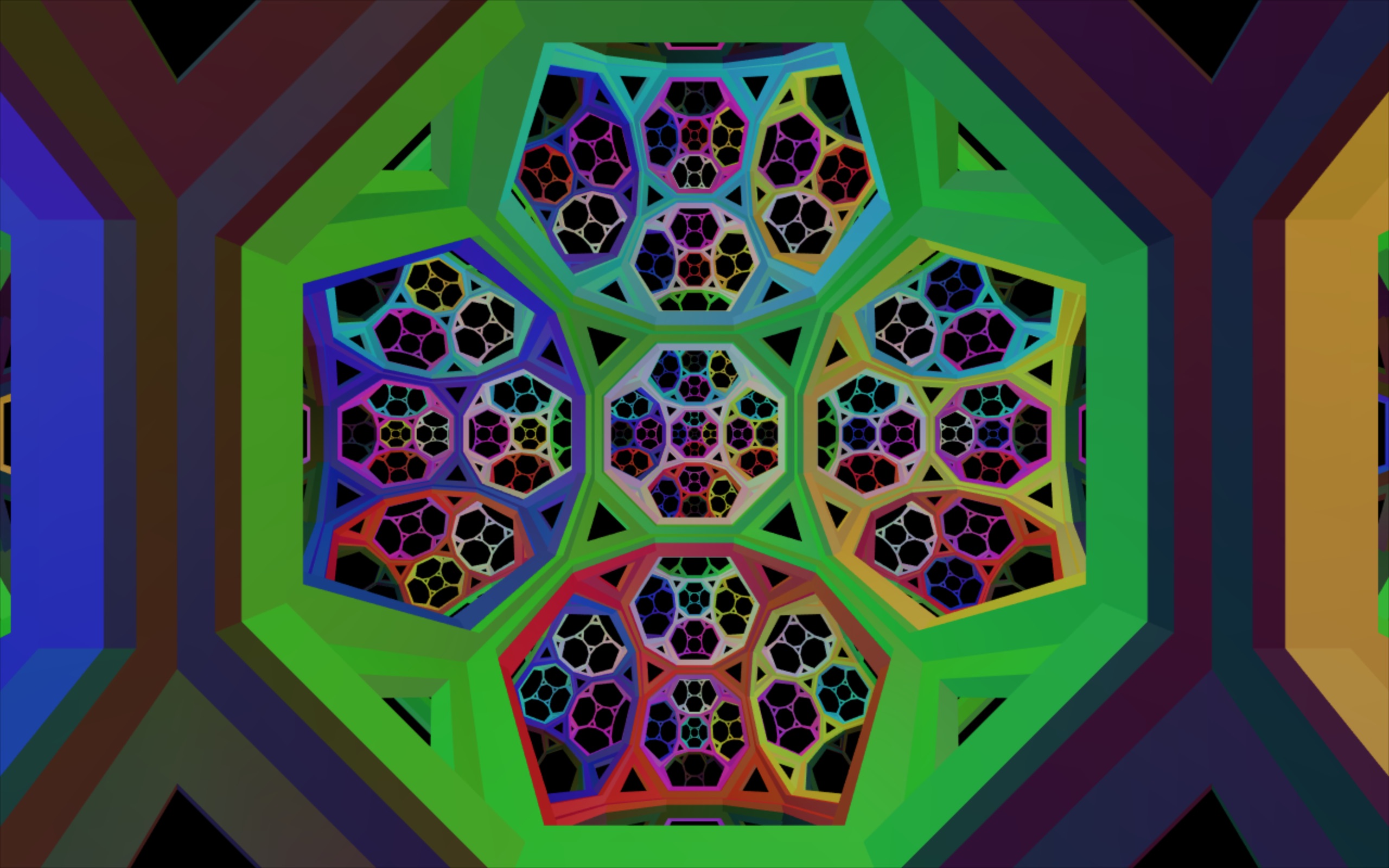}
\label{H3_parallel_transport_1}
}
\quad
\subfloat[$\HH^2\times \EE$ view 1: along $\EE$.]
{
\includegraphics[width=0.22\textwidth]{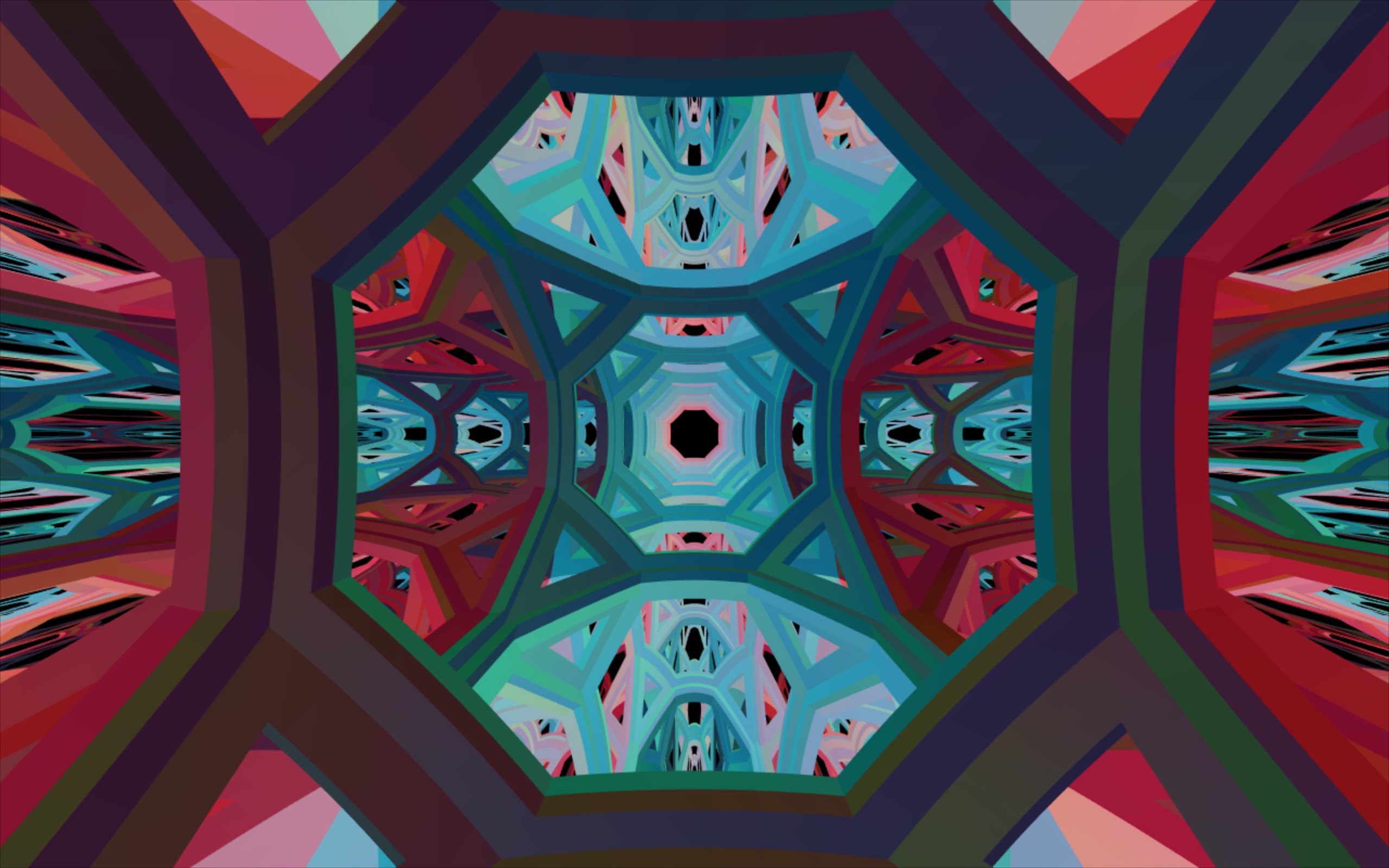}
\label{H2xR_parallel_transport1_1}
}
\quad
\subfloat[$\HH^2\times \EE$ view 2: $45^\circ$ between the $\EE$ and $\HH^2$ directions.]
{
\includegraphics[width=0.22\textwidth]{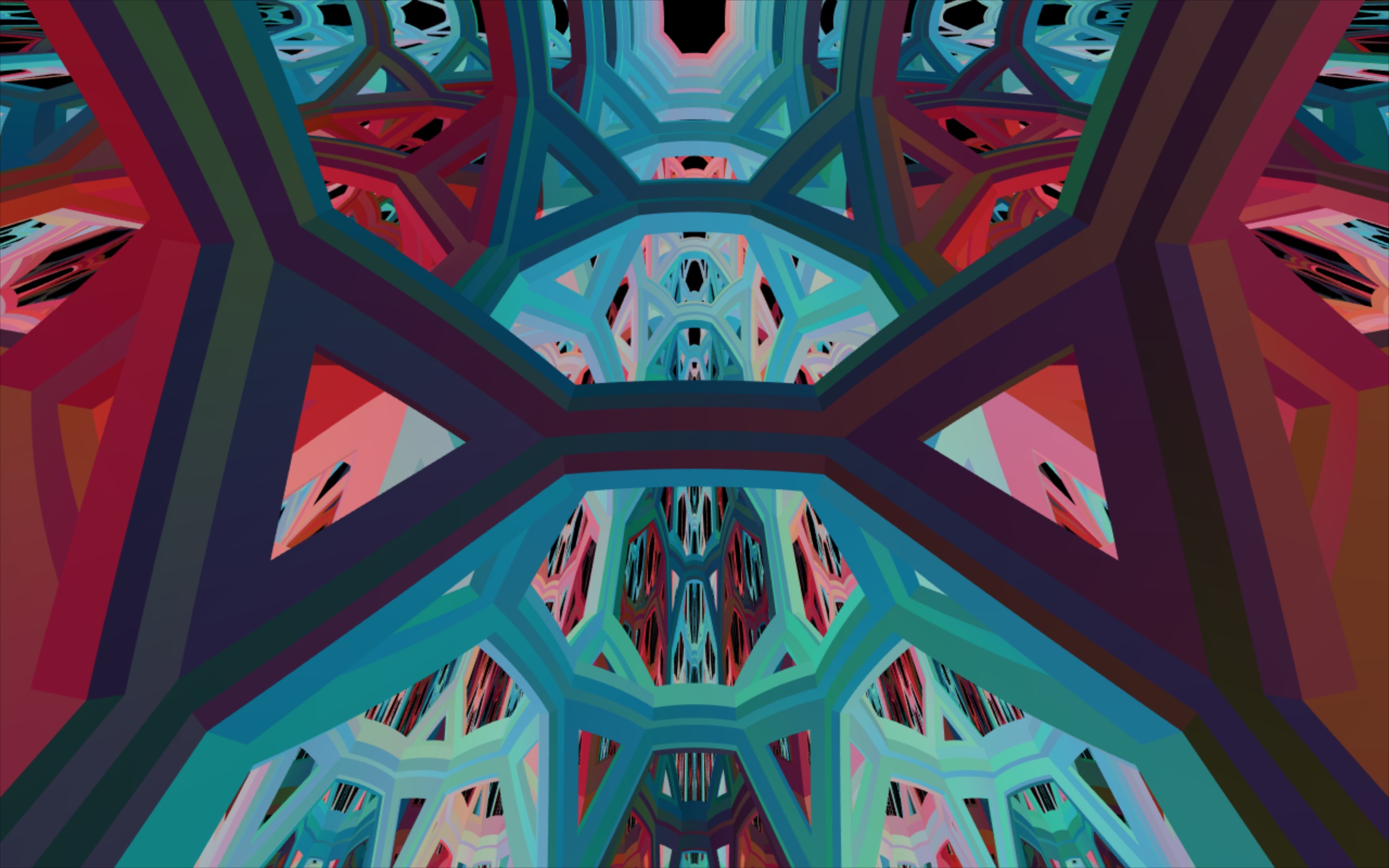}
\label{H2xR_parallel_transport2a_1}
}
\quad
\subfloat[$\HH^2\times \EE$ view 3: along $\HH^2$.]
{
\includegraphics[width=0.22\textwidth]{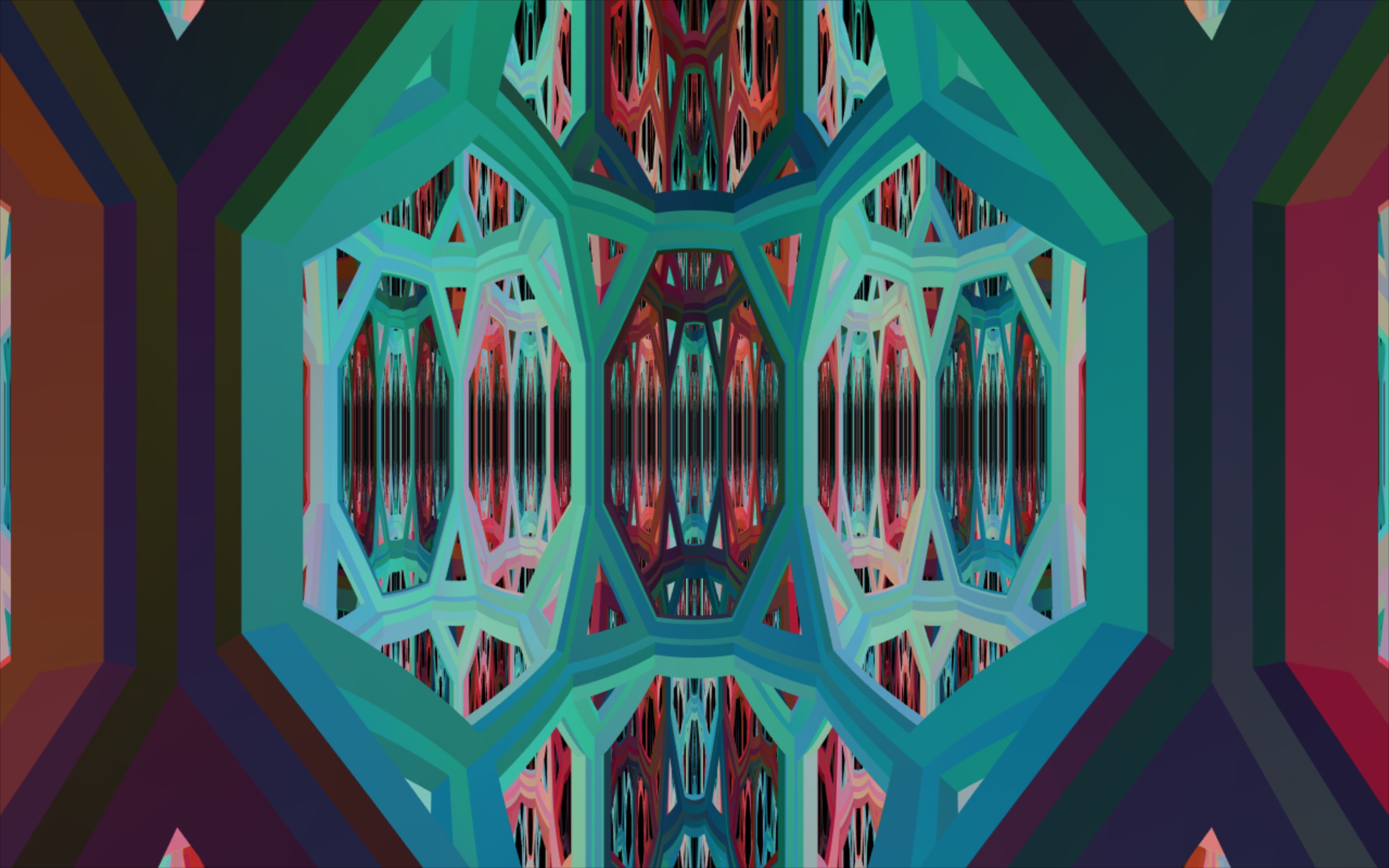}
\label{H2xR_parallel_transport3_1}
}

\vspace{-7pt}
\subfloat[$\HH^3$ after moving right 0.5.]
{
\includegraphics[width=0.22\textwidth]{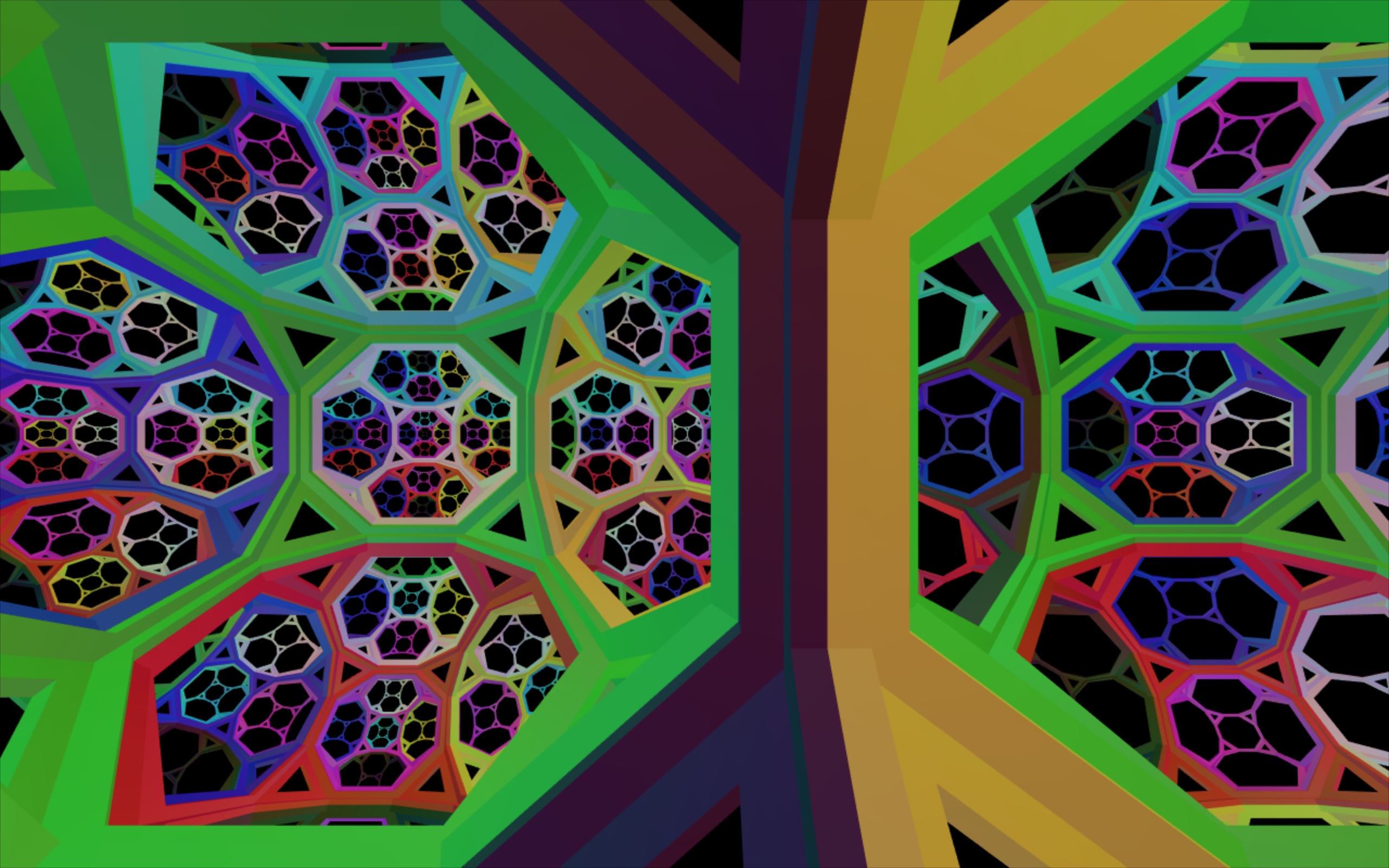}
\label{H3_parallel_transport_2}
}
\quad
\subfloat[$\HH^2\times \EE$ view 1 after moving right 0.5.]
{
\includegraphics[width=0.22\textwidth]{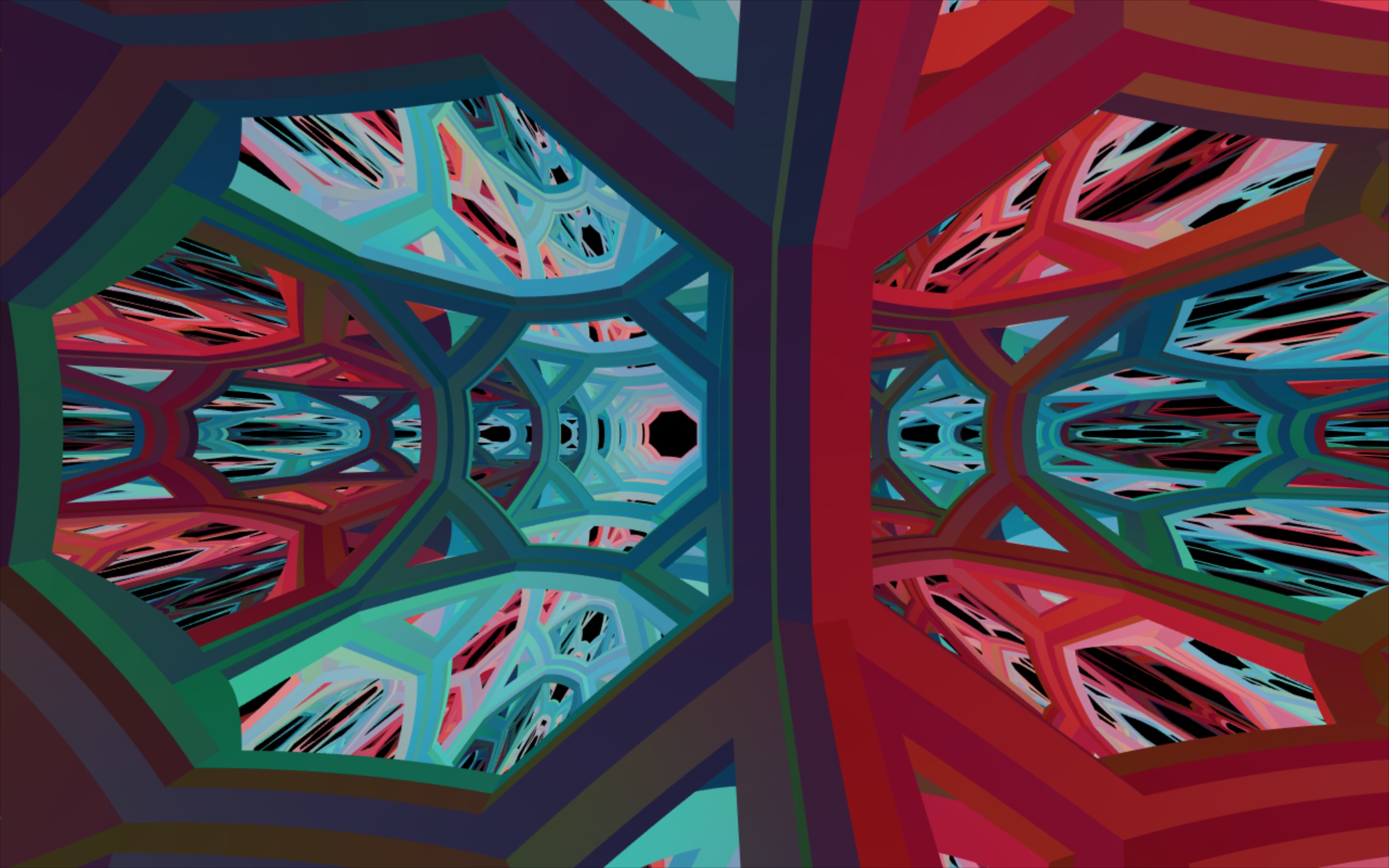}
\label{H2xR_parallel_transport1_2}
}
\quad
\subfloat[$\HH^2\times \EE$ view 2 after moving right 0.5.]
{
\includegraphics[width=0.22\textwidth]{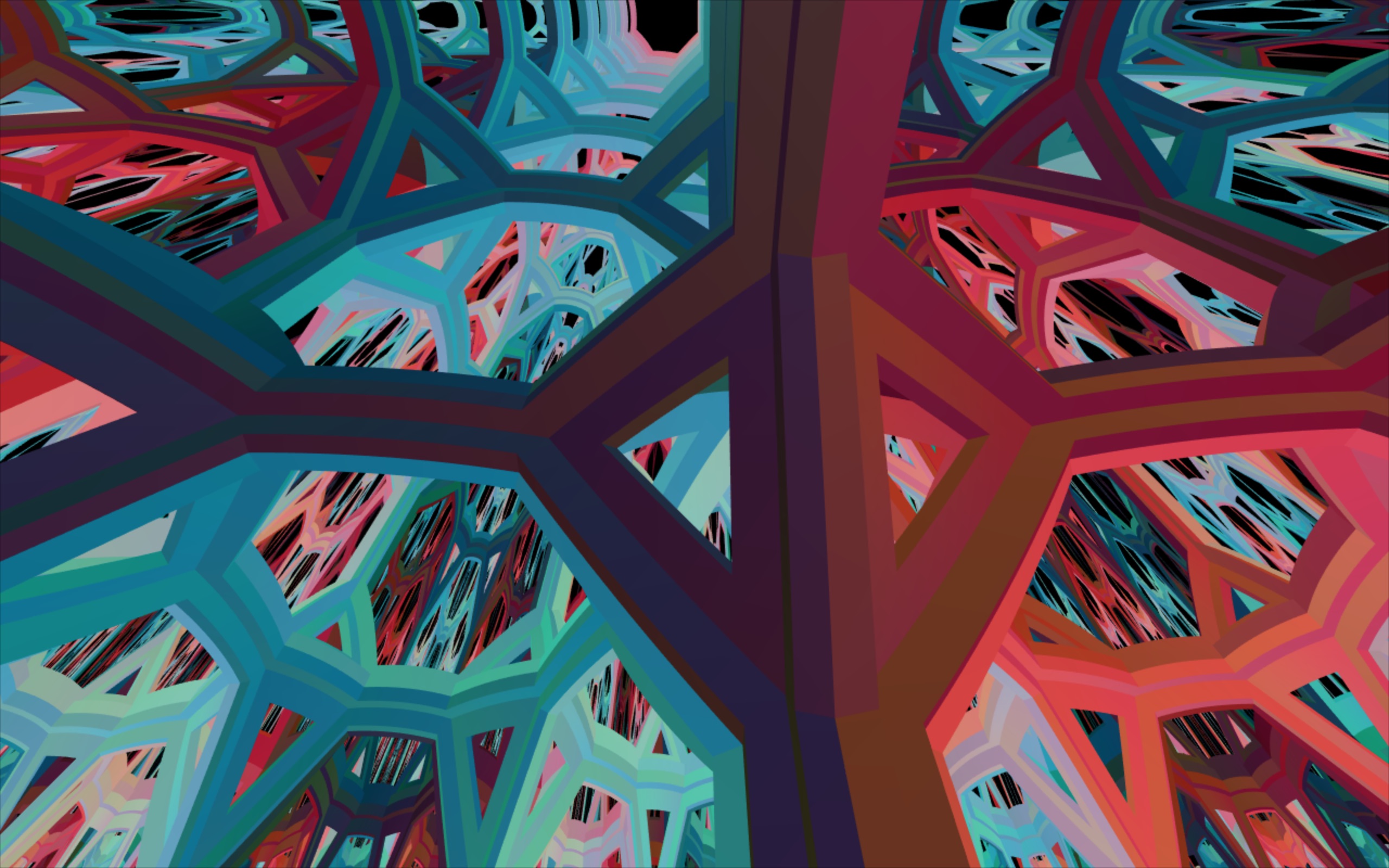}
\label{H2xR_parallel_transport2a_2}
}
\quad
\subfloat[$\HH^2\times \EE$ view 3 after moving right 0.5.]
{
\includegraphics[width=0.22\textwidth]{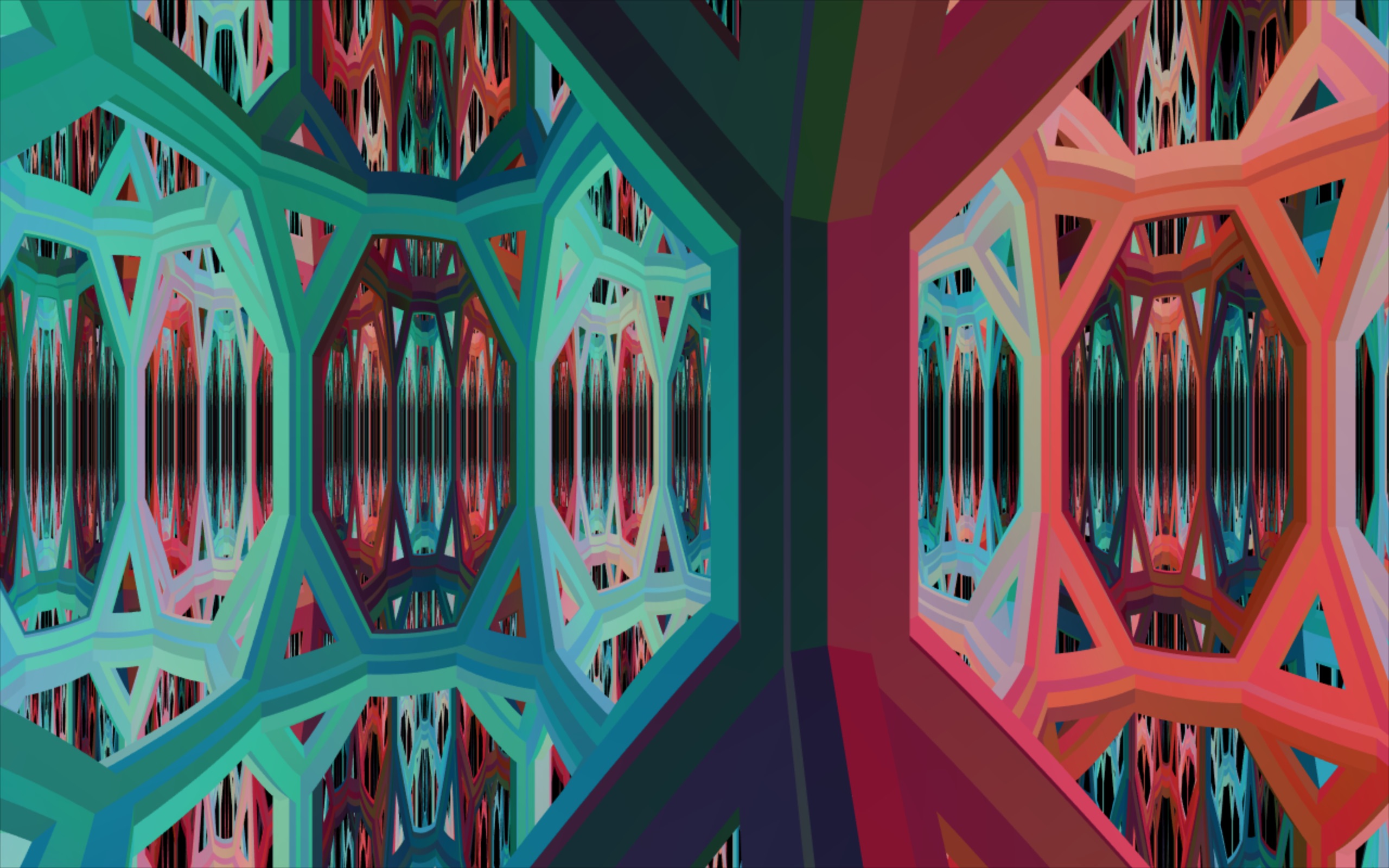}
\label{H2xR_parallel_transport3_2}
}

\vspace{-7pt}
\subfloat[$\HH^3$ after moving up 0.5.]
{
\includegraphics[width=0.22\textwidth]{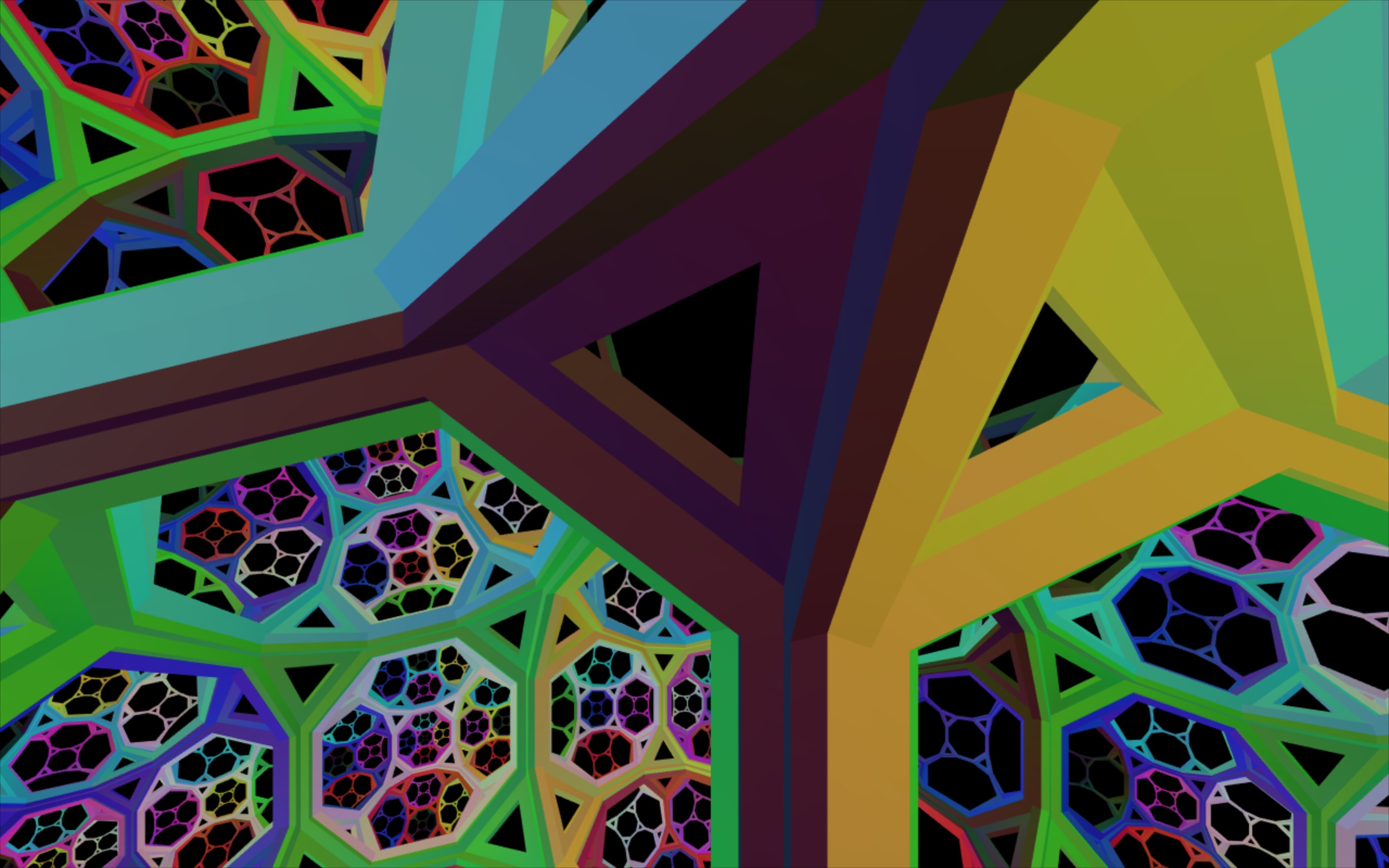}
\label{H3_parallel_transport_3}
}
\quad
\subfloat[$\HH^2\times \EE$ view 1 after moving up 0.5.]
{
\includegraphics[width=0.22\textwidth]{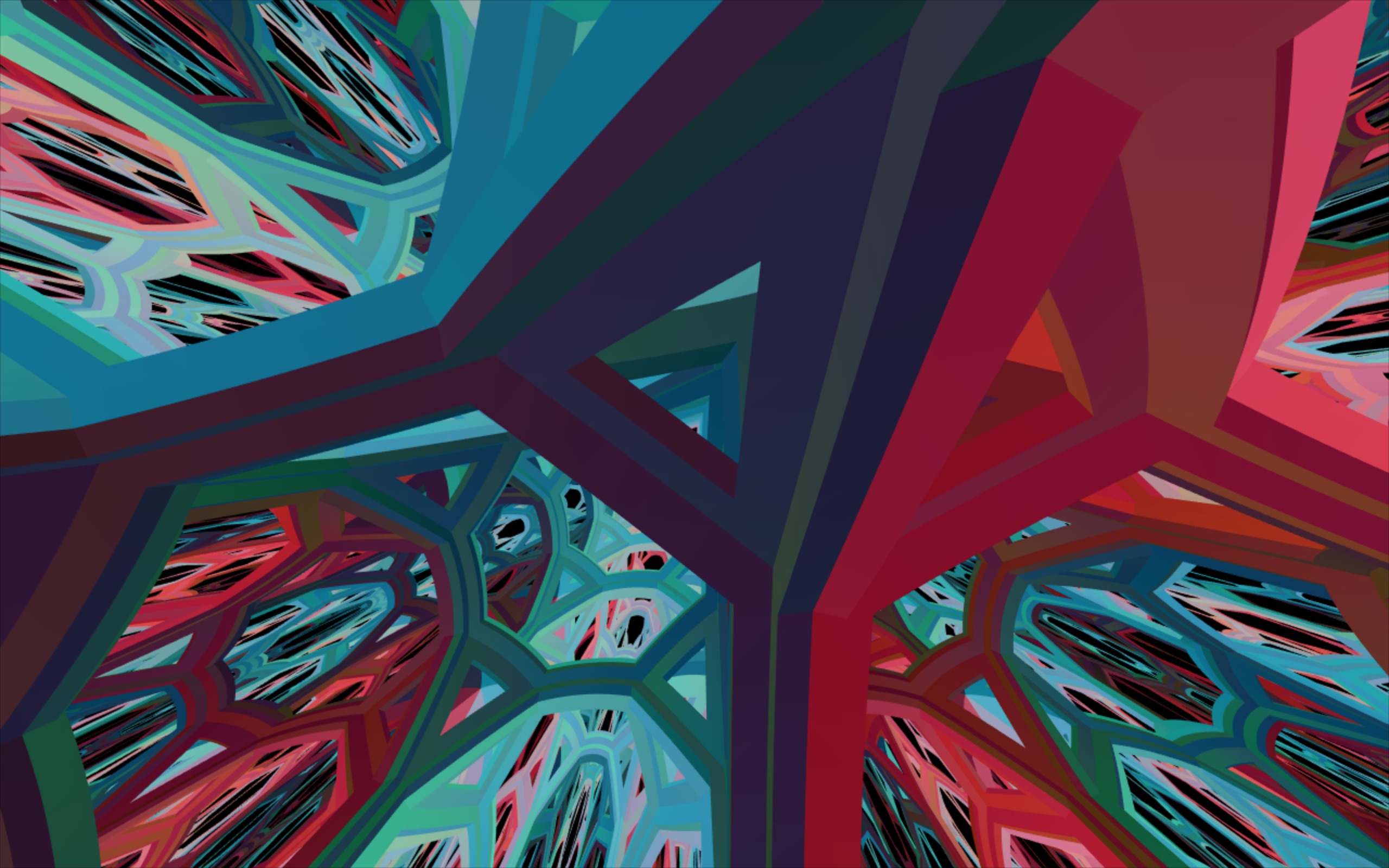}
\label{H2xR_parallel_transport1_3}
}
\quad
\subfloat[$\HH^2\times \EE$ view 2 after moving up 0.5.]
{
\includegraphics[width=0.22\textwidth]{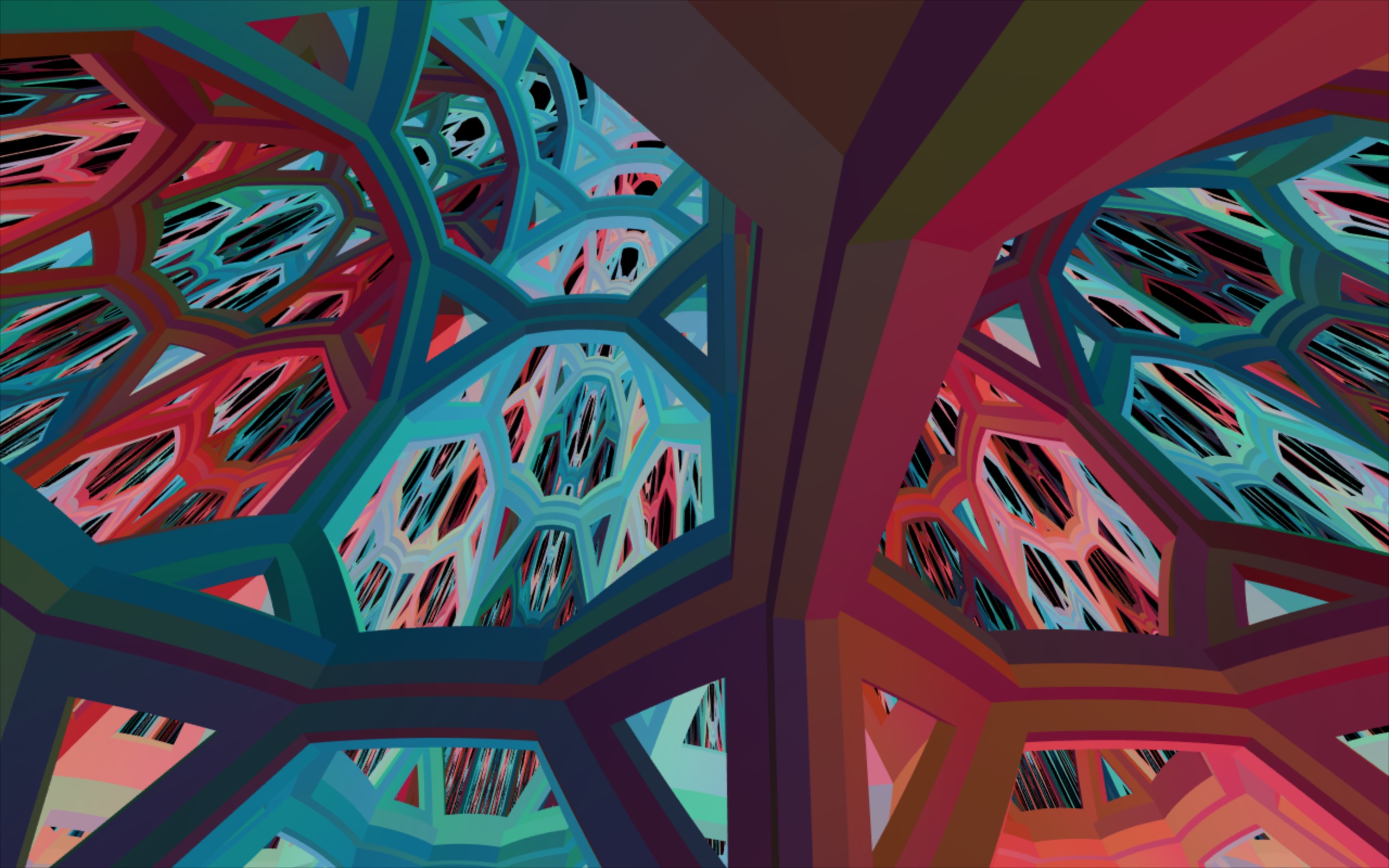}
\label{H2xR_parallel_transport2a_3}
}
\quad
\subfloat[$\HH^2\times \EE$ view 3 after moving up 0.5.]
{
\includegraphics[width=0.22\textwidth]{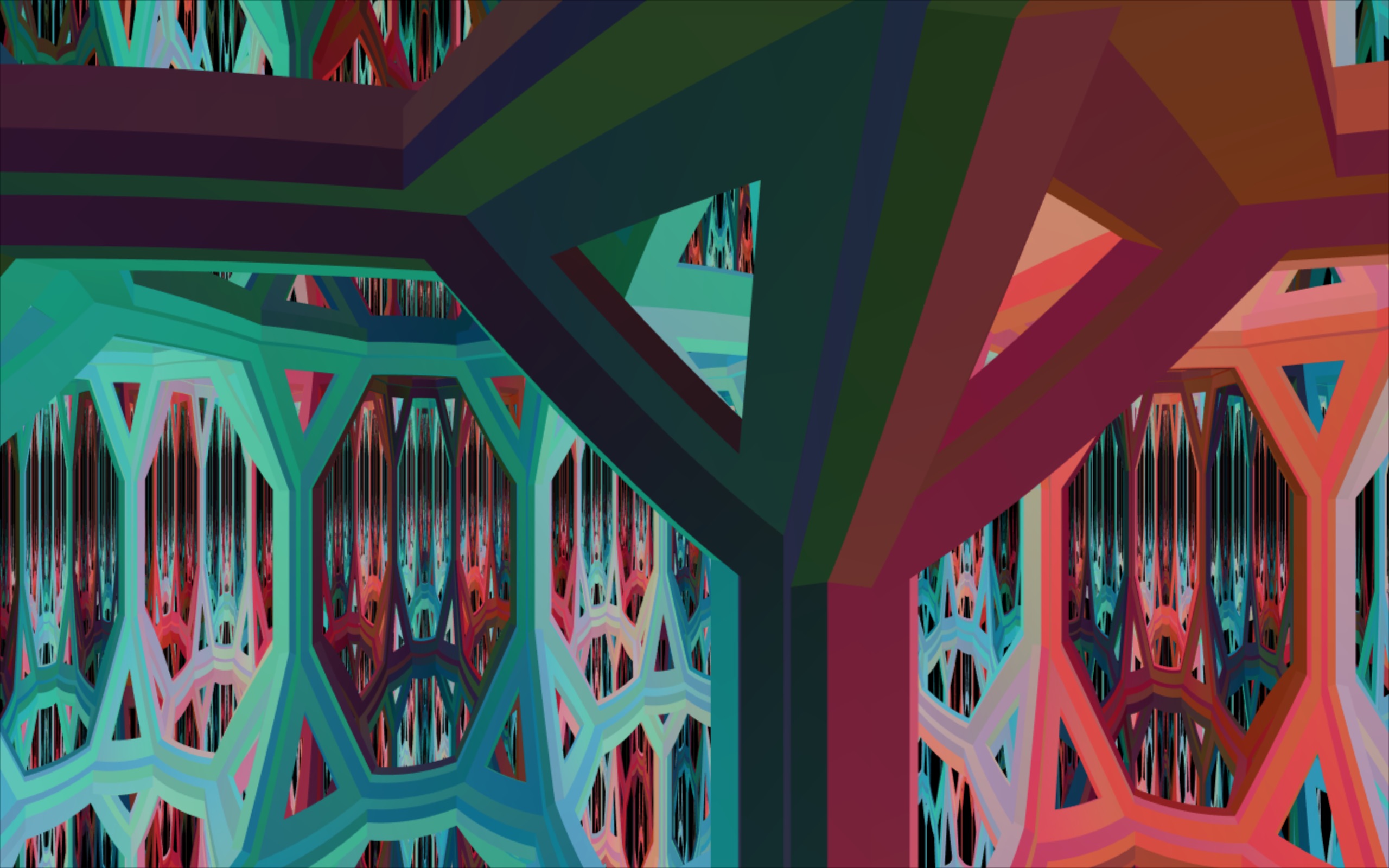}
\label{H2xR_parallel_transport3_3}
}

\vspace{-7pt}
\subfloat[$\HH^3$ after moving left 0.5.]
{
\includegraphics[width=0.22\textwidth]{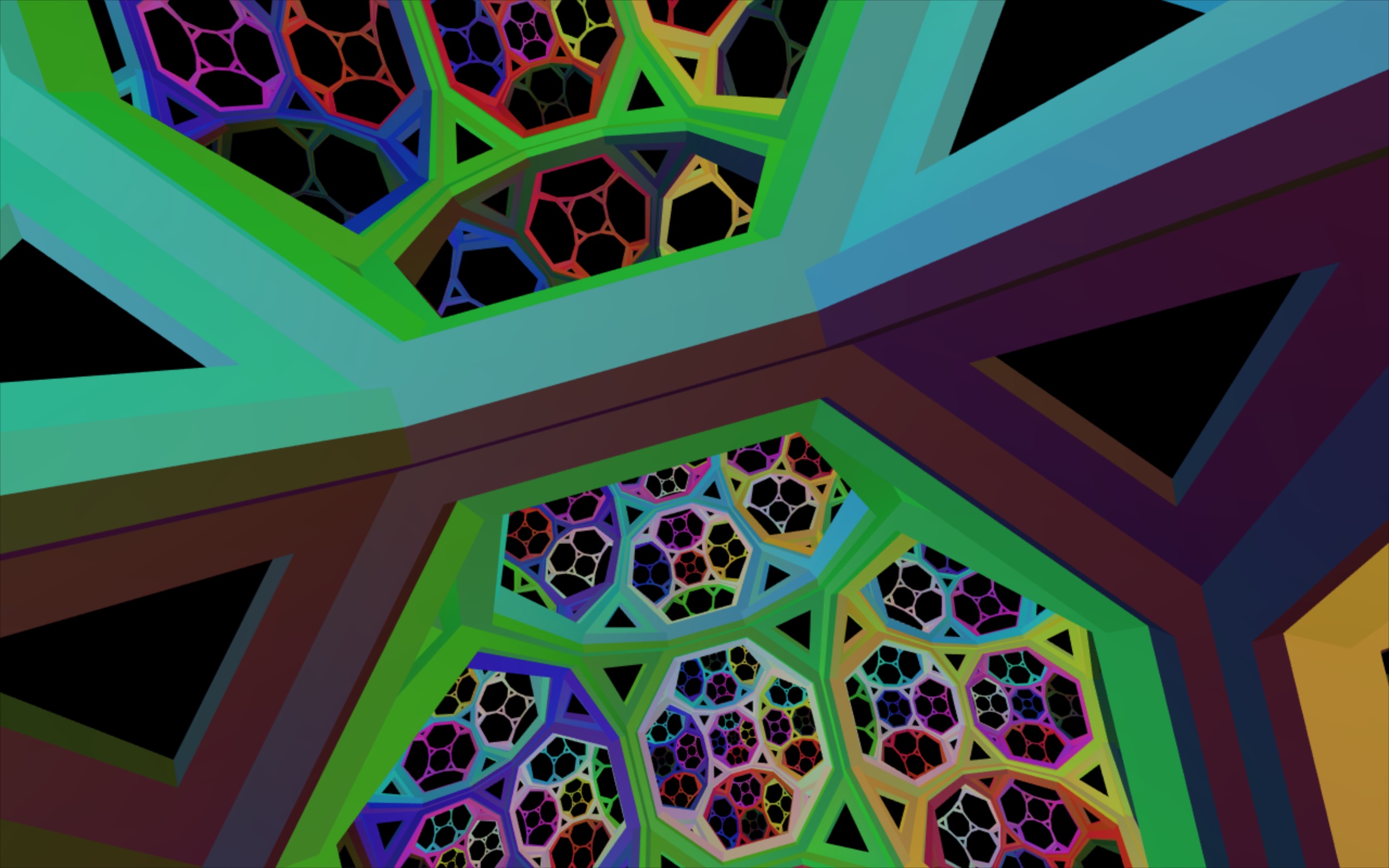}
\label{H3_parallel_transport_4}
}
\quad
\subfloat[$\HH^2\times \EE$ view 1 after moving left 0.5.]
{
\includegraphics[width=0.22\textwidth]{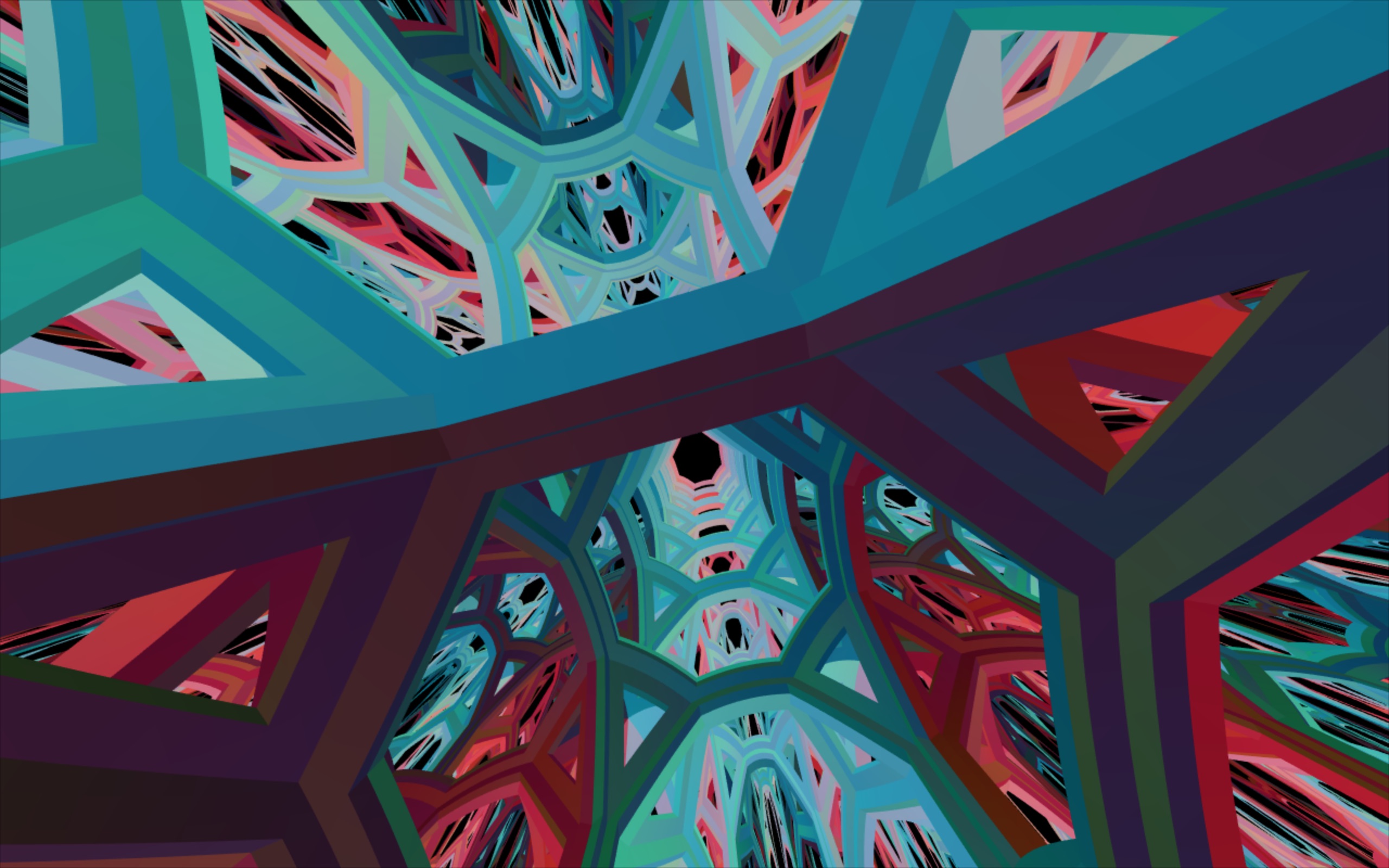}
\label{H2xR_parallel_transport1_4}
}
\quad
\subfloat[$\HH^2\times \EE$ view 2 after moving left 0.5.]
{
\includegraphics[width=0.22\textwidth]{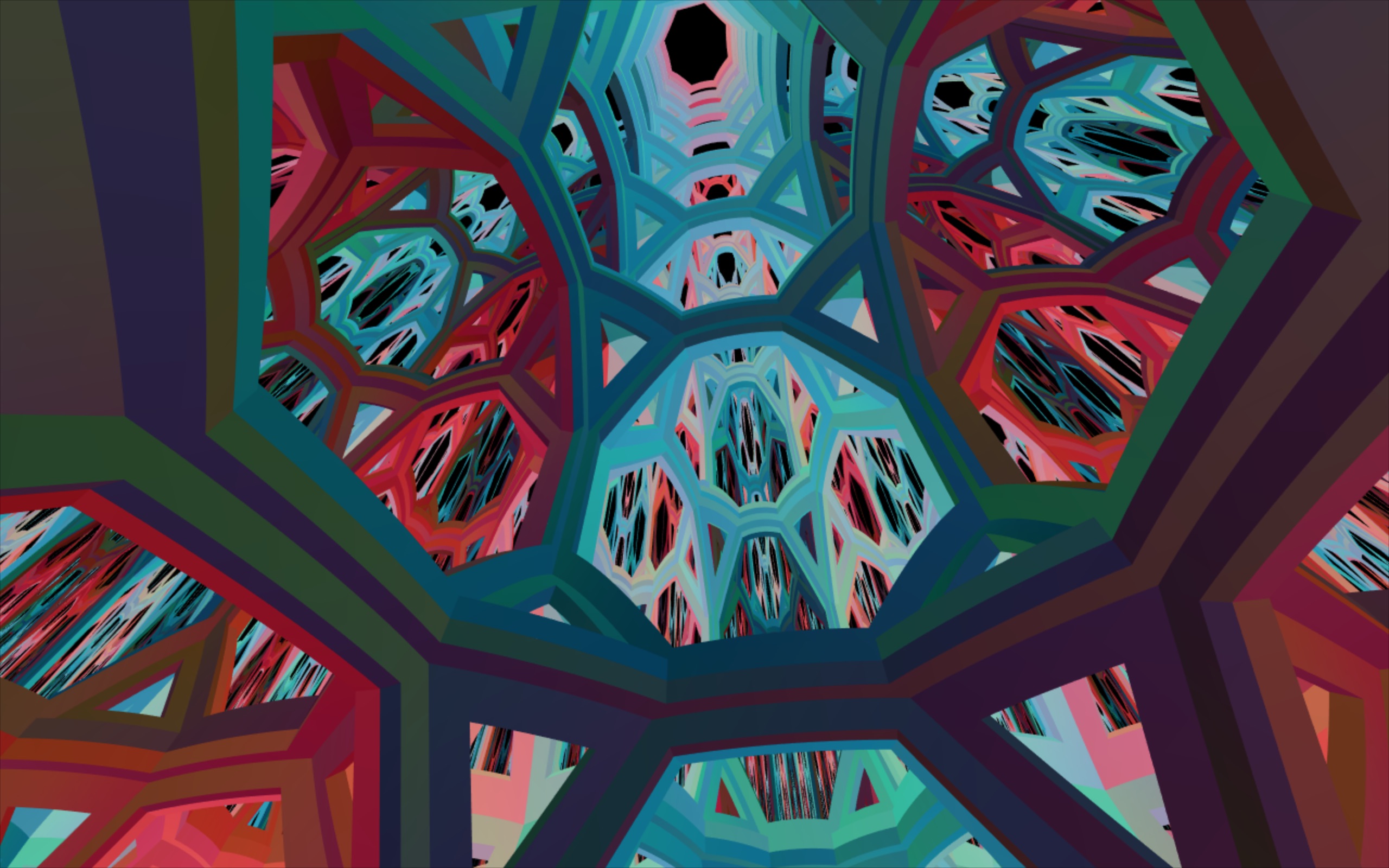}
\label{H2xR_parallel_transport2a_4}
}
\quad
\subfloat[$\HH^2\times \EE$ view 3 after moving left 0.5.]
{
\includegraphics[width=0.22\textwidth]{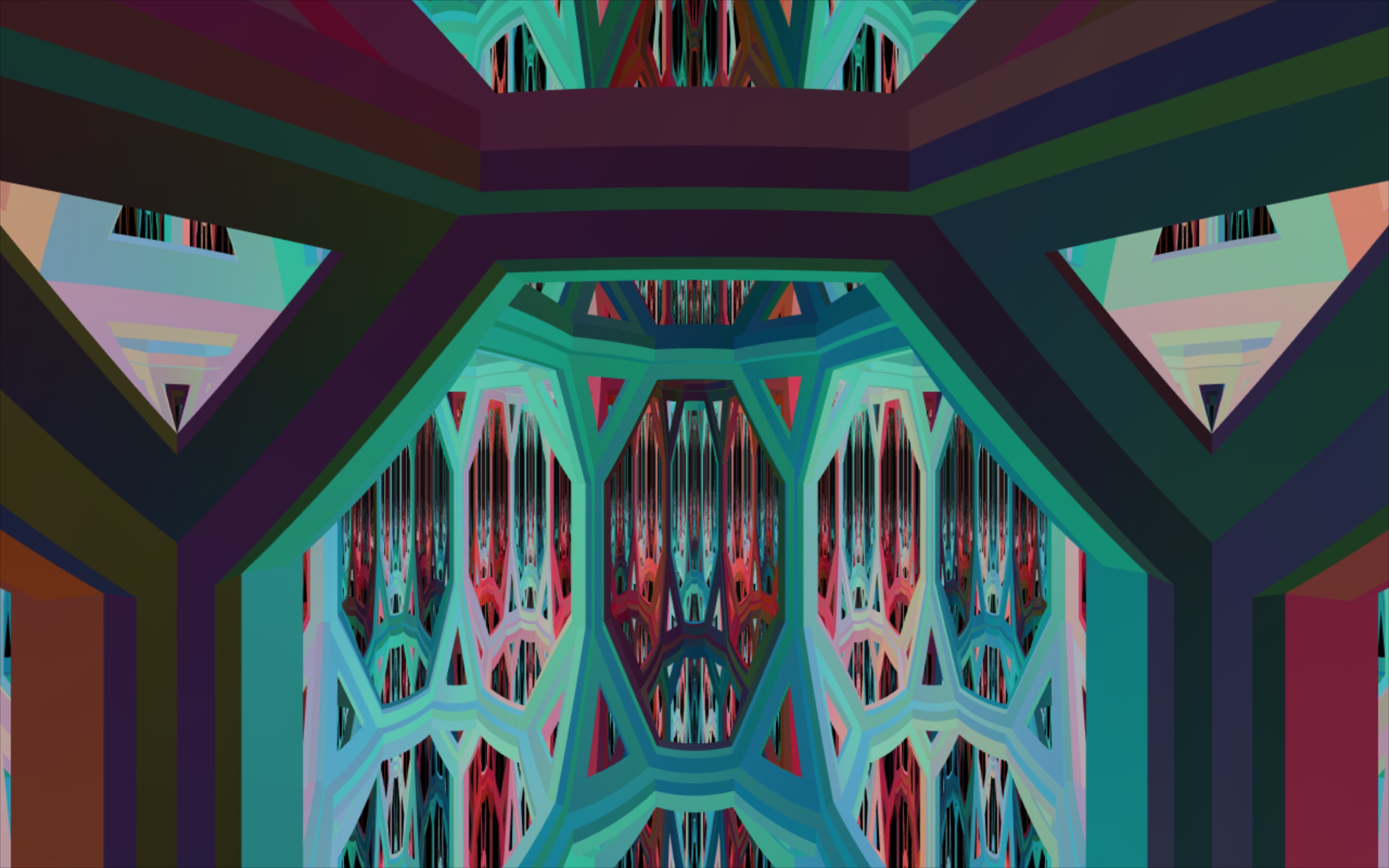}
\label{H2xR_parallel_transport3_4}
}

\vspace{-7pt}
\subfloat[$\HH^3$ after moving down 0.5.]
{
\includegraphics[width=0.22\textwidth]{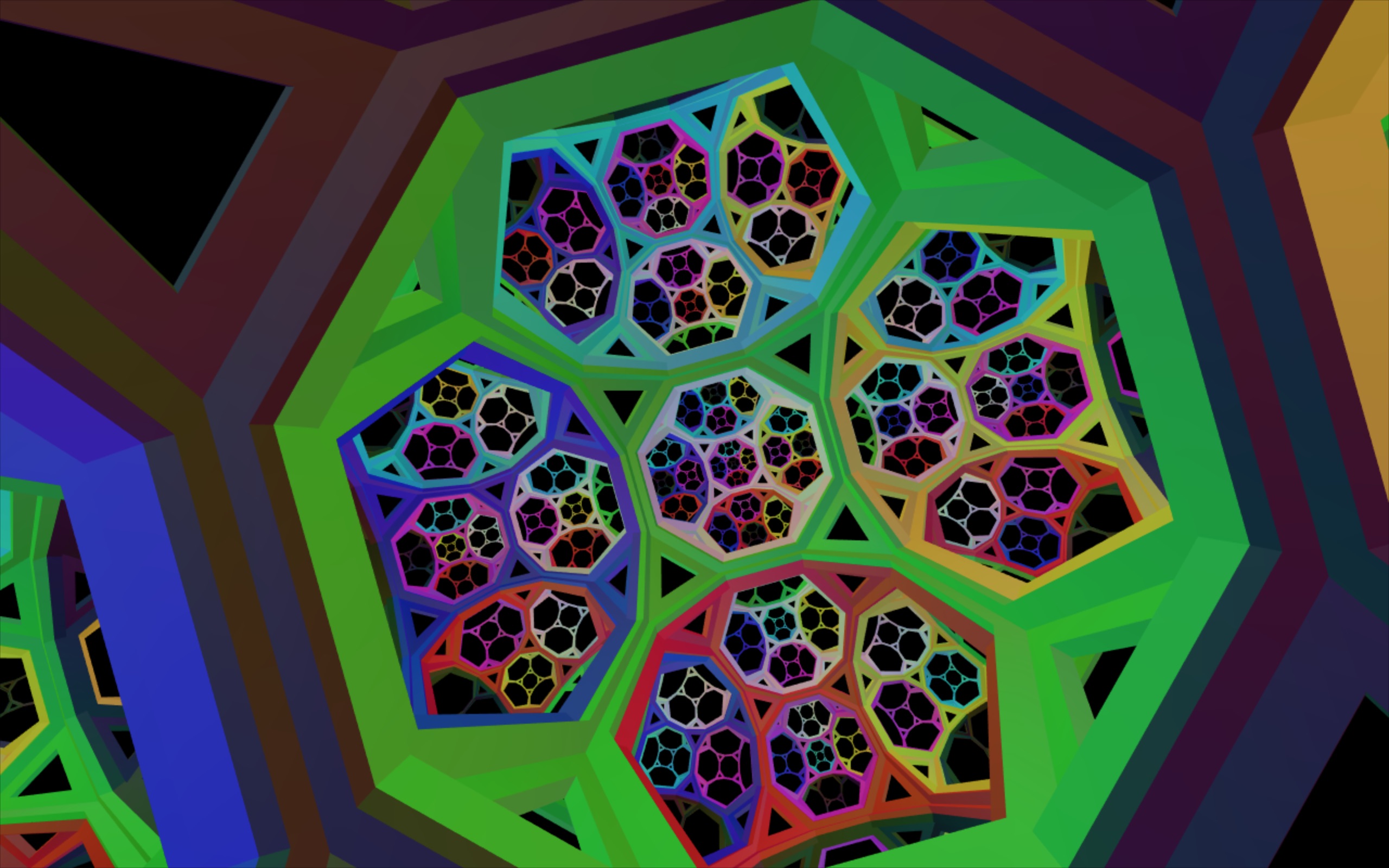}
\label{H3_parallel_transport_5}
}
\quad
\subfloat[$\HH^2\times \EE$ view 1 after moving down 0.5.]
{
\includegraphics[width=0.22\textwidth]{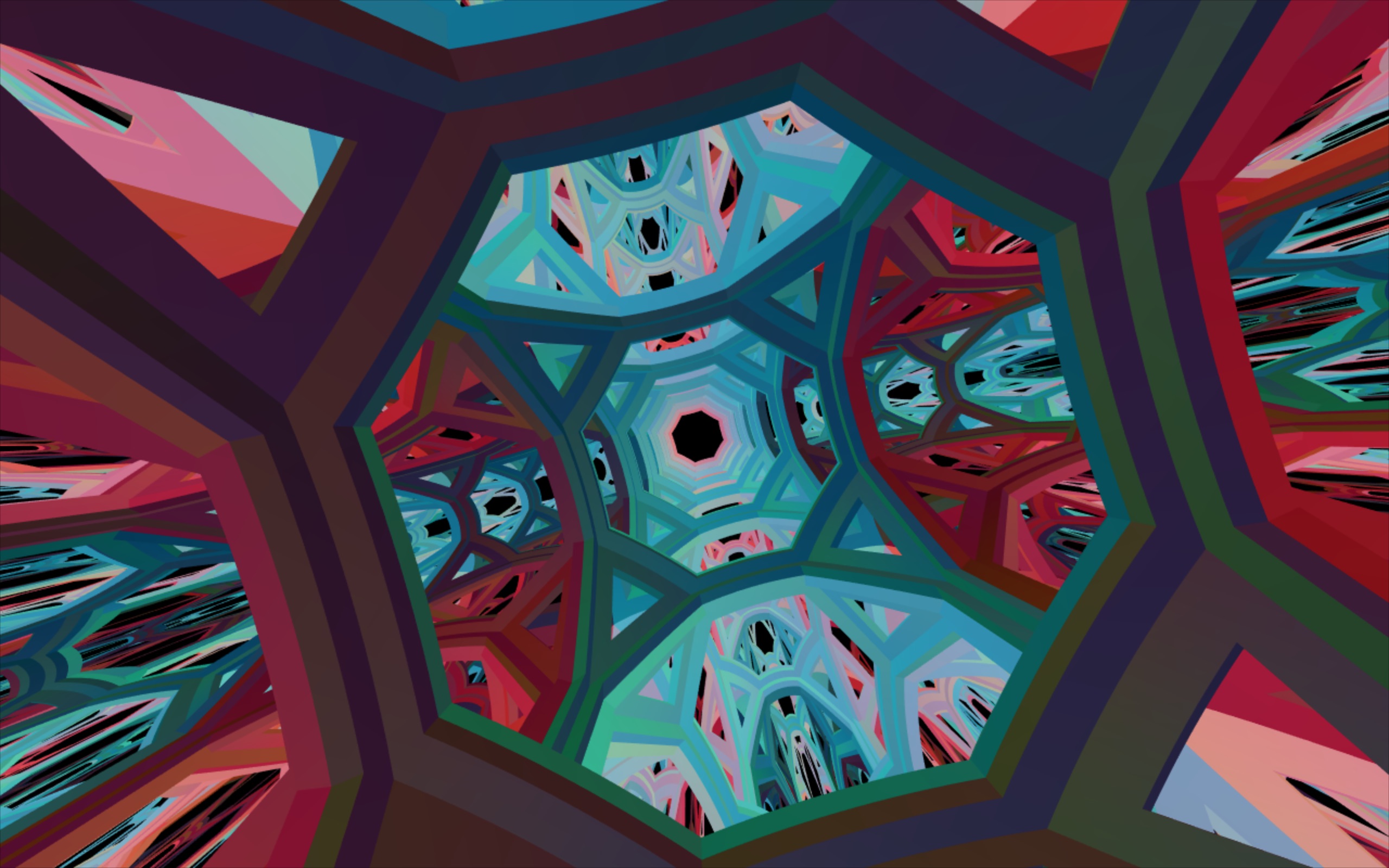}
\label{H2xR_parallel_transport1_5}
}
\quad
\subfloat[$\HH^2\times \EE$ view 2 after moving down 0.5.]
{
\includegraphics[width=0.22\textwidth]{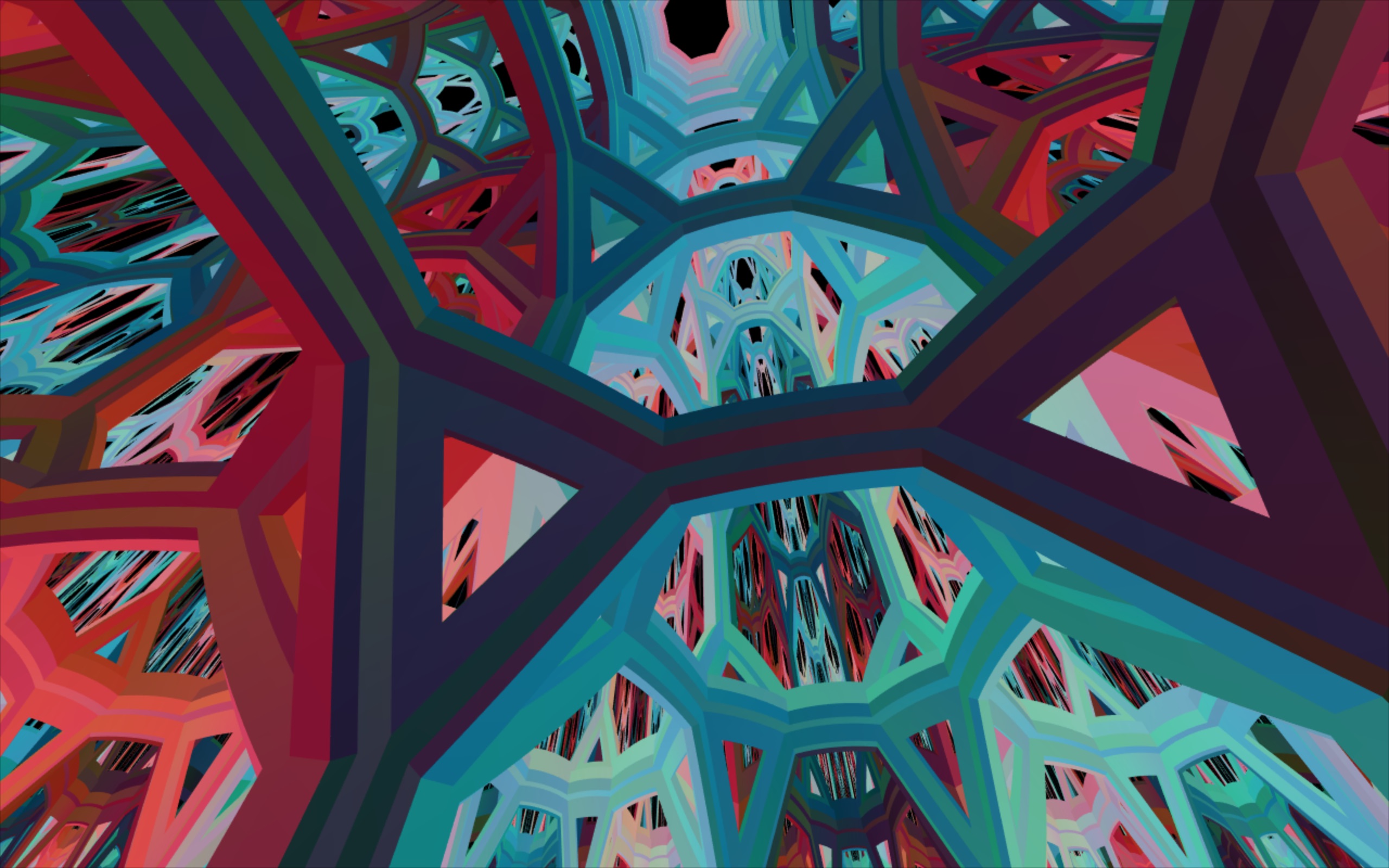}
\label{H2xR_parallel_transport2a_5}
}
\quad
\subfloat[$\HH^2\times \EE$ view 3 after moving down 0.5.]
{
\includegraphics[width=0.22\textwidth]{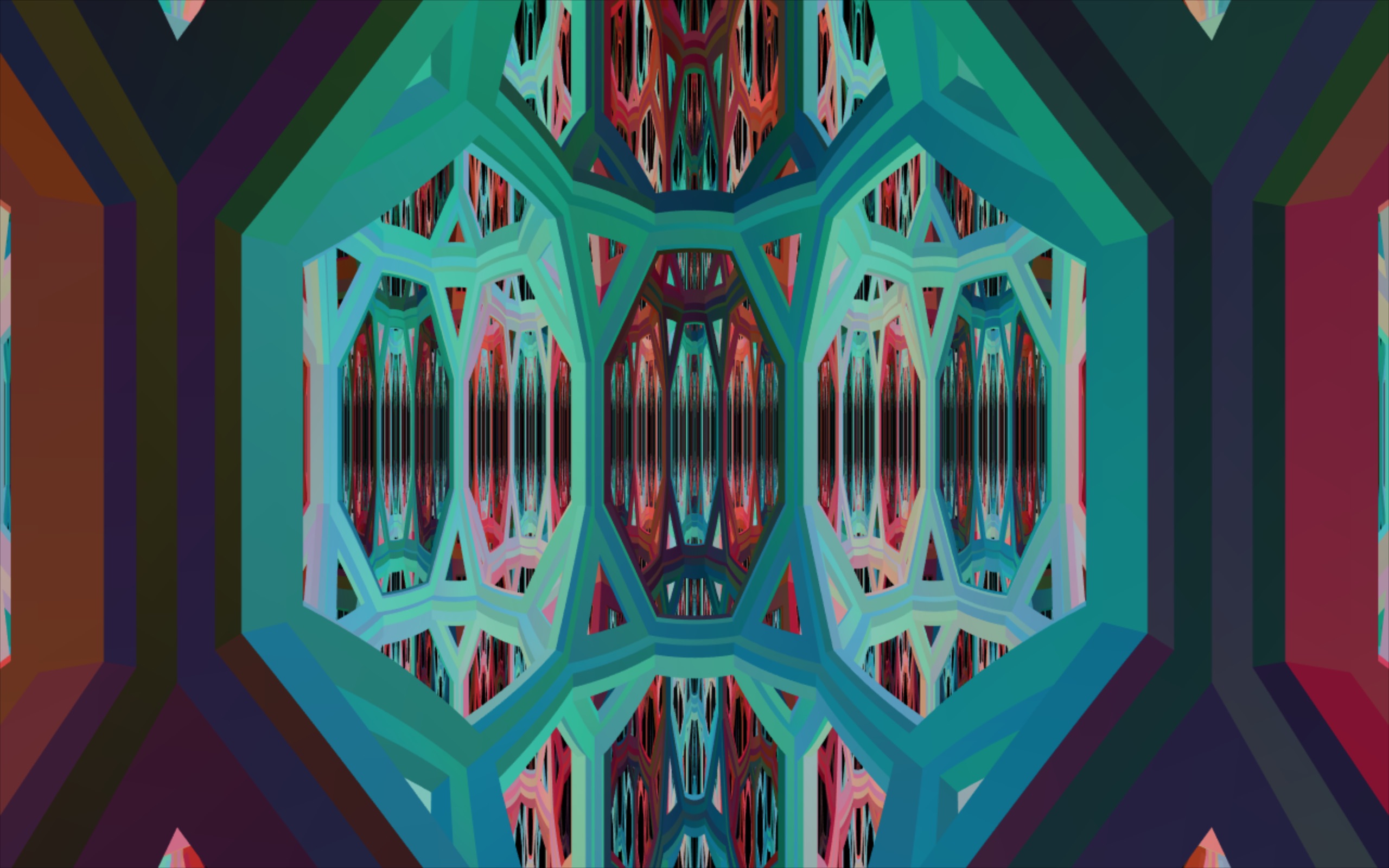}
\label{H2xR_parallel_transport3_5}
}
\caption{Parallel transport rotates reference frames in curved space. Rotations in the $\HH^2$ plane (column 2) have the same effects as looking in any direction in $\HH^3$ (column 1). However other directions behave differently.}
\label{Fig:parallel_transport}
\vspace{-5pt}
\end{figure}

Eventually, we plan to put recognisable, human scale objects in our simulation of $\HH^2 \times \EE$, for people to interact with. For now, we pattern it with cubes, as we did for $\HH^3$~\cite{our_paper_h3}. In each horizontal $\HH^2$ slice, we draw the $\{4,6\}$ Schl\"afli symbol tiling: squares (4 sides), with 6 meeting at each vertex. These horizontal squares are cross-sections of our cubes. At regular intervals in the $\EE$ direction, one cube ends, and the next begins. Because this is a product space, we have freedom to choose the height of the cubes in $\EE$ relative to the side length of the $\{4,6\}$ tiling of $\HH^2$. We choose this height so that when viewed from its centre, a cube in $\HH^2\times\EE$ has its vertices coincide with the vertices of a euclidean cube. The inverse of the exponential map associates each of these vertices to a point in the tangent space that is a distance $\arcsinh(\sqrt{2})/\sqrt{2}$ away from the centre. 
We take the vertices in the meshes that make up the cubes and their decorations in $\EE^3$ and move them upstairs into 
(hyperboloid)$\times\EE$ via the map from the Klein model to the hyperboloid model according to the coordinate transformation $(x,y,z)\rightarrow(x/\sqrt{1-x^2-y^2},y/\sqrt{1-x^2-y^2},z,1/\sqrt{1-x^2-y^2})$. Once in the hyperboloid model, 
transformations are handled through the exponential map. 
Unlike $\HH^3$ which is isotropic, the anisotropy of $\HxE$ causes the 
straight lines bounding the faces of the cube in $\EE^3$
to become curved as you look along a 
diagonal direction in $\HxE$.

In Figure \ref{h2xr_46}, we colour each layer of cubes in essentially the same way, coming from a colouring of the $\{4,6\}$ tiling. Subsequent layers have their colours ``rotated'' slightly, around the circle of colours from black to red to white to cyan. As in our colouring of the $\{4,3,6\}$ honeycomb, the colouring on each layer comes from lifting a colouring of a closed manifold, in this case the genus two surface. See Figure \ref{46_colouring}.

\section{Revenge of curved spaces in virtual reality}

As with $\HH^3$~\cite{our_paper_h3}, our experience in $\HxE$ differs viscerally from in $\EE^3$. One of the most striking differences between how objects appear in $\HxE$ versus how objects appear in $\EE^3$ and $\HH^3$ stems from its lack of isotropy. As we move forwards, the aspect ratio of objects in front of us changes. Figure \ref{Fig:aspect_ratio} shows a sequence of views, moving forwards from one cube into the next. The octagon in front of us starts out tall and thin, and becomes wider as we approach. The octagons further into the distance are even thinner. The reason for this behaviour in $\HxE$ is that an object's height goes down linearly relative to its distance from us (just as it does in $\EE^3$) 
but its width scales exponentially. This doesn't happen in $\HH^3$ -- both width and height scale exponentially as $\HH^3$ is isotropic, so the aspect ratio of objects does not change.

An interesting feature of $\HxE$ being a product space is that parallel transport affects directions differently. In $\HH^3$ when you traverse a closed loop in any plane, you experience the world as having rotated. This is because the vector pointing ``up" gets parallel transported around the loop and comes back as rotated~\cite{our_paper_h3}. In $\HxE$, you experience this phenomenon when looking along the $\EE$ direction -- that is to say, the closed loop lies in the $\HH^2$ plane. See Figure \ref{Fig:parallel_transport}. The frame is transported in the same it was in any direction in $\HH^3$. However, if you look directly along the $\HH^2$ plane and traverse a closed loop, first moving to the right for distance $d$ along the $\HH^2$ direction, then up along $\EE$ for distance $d$, then to the left along $\HH^2$ and finally returning to the starting point by travelling down along $\EE$, you will notice that the world has not rotated. This is because pure isometries in $\HH^2$ and $\EE$ commute as $\HxE$ is a cartesian product.

\subsection{Parallax}

\begin{wrapfigure}[12]{l}{0.5\textwidth}
\vspace{-5pt}
\centering
\subfloat[Geodesics (light rays) from two eyes looking straight ahead would diverge in $\HH^2$.]
{
\includegraphics[width=0.225\textwidth]{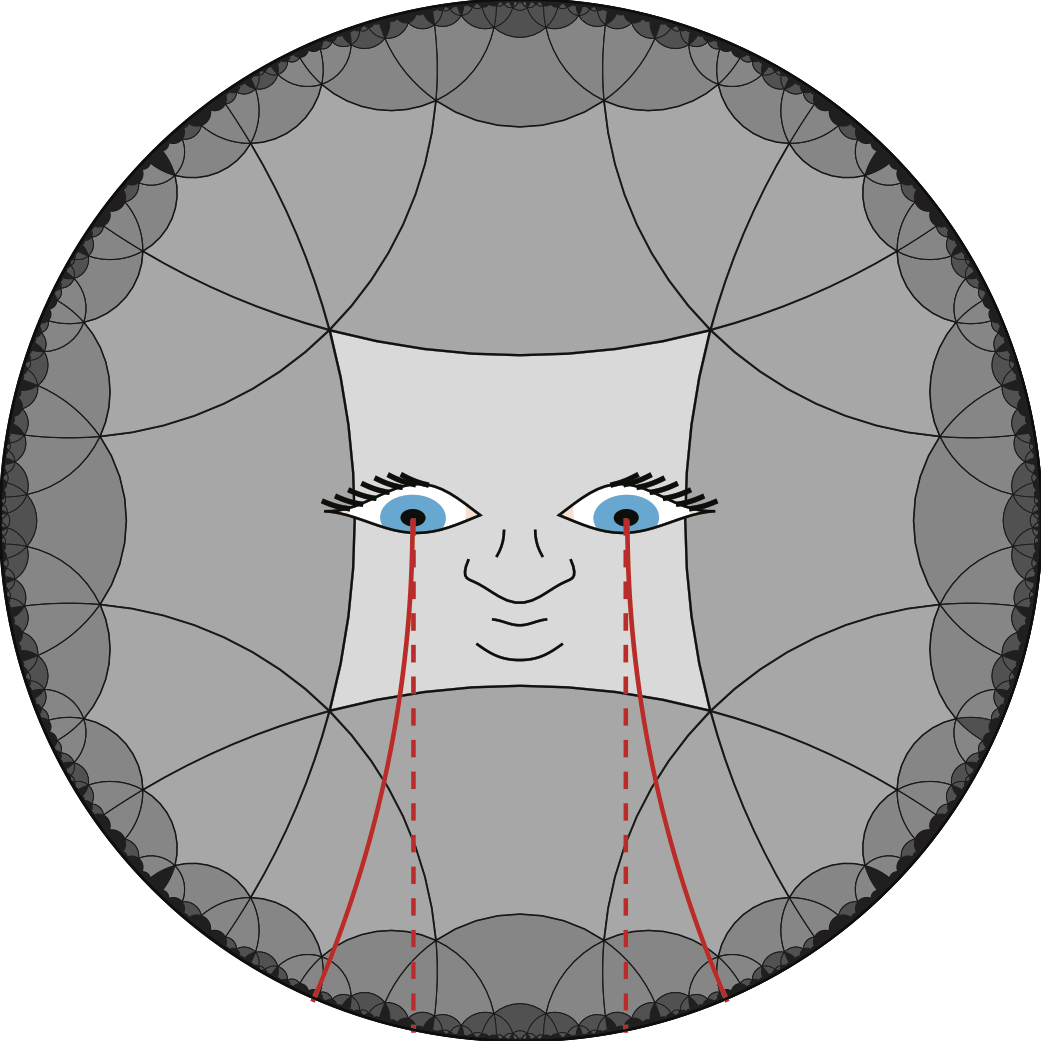}
\label{parallax1_1}
}
\quad
\subfloat[Eyes must point inward to see an object that is infinitely far away in $\HH^2$.]
{
\includegraphics[width=0.225\textwidth]{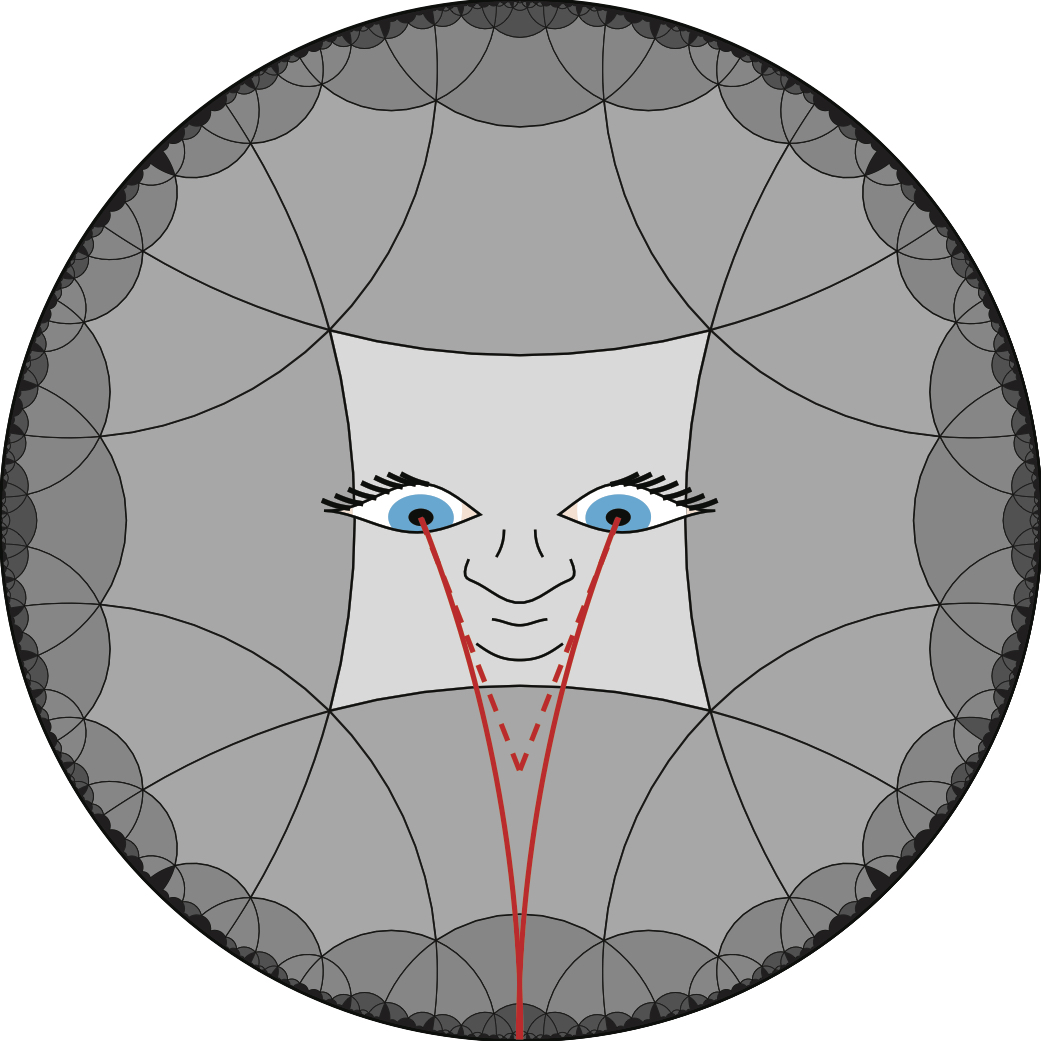}
\label{parallax_2}
}
\caption{Parallel transport is responsible for the phenomena associated with parallax in $\HH^2$. See also~\cite{weeks_real-time_rendering}.}
\label{Fig:creepy_eyes}
\end{wrapfigure}

Humans judge the distance of objects from them using a number of different techniques. One such technique is our use of \emph{stereopsis}, the perception of depth inferred from differences in the visual input coming into our two eyes. This is an example of a use of \emph{parallax}, the difference between the apparent position of an object viewed from different viewpoints. An object that is directly in front of us and close to us, will appear on the right side of the field of view of one's left eye, and on the left side of the field of view of one's right eye. An object that is directly in front of us but very far away from us will appear in the center of the field of views from each eye, assuming that our eyes are both pointing directly forward. See Figure \ref{Fig:creepy_eyes}.

In particular, an inhabitant of euclidean space expects to see an object that is infinitely far away along rays of light that enter their two eyes along parallel geodesics. In negatively curved spaces (e.g. $\HH^3$ and $\HH^2$)  however, geodesics that enter one's eyes along parallel directions diverge as we follow the light rays backwards -- there can be no single object that those light rays came from. Instead, a euclidean visitor to hyperbolic space has to turn their eyes inwards to see an object that is infinitely far away. The effect is that 
all objects appear to be relatively close by. In general, objects in hyperbolic space are further away than they appear.

\section{Future directions: Moving objects in the space}

Moving objects in $\HxE$ presents several interesting challenges. In addition to the headset, the HTC Vive can also track the position and orientation of hand-held controllers. In most applications that use the controller, a virtual version of the controller is visible in the virtual space. This helps greatly with the user's sense of embodiment in the space, since they can see the positions of their hands. The controller also has various buttons and triggers for other forms of interaction, for example grabbing on to a virtual object near to the virtual position of the controller, allowing the object to be moved in space. We plan to add this kind of interactive element to our simulations. 

An obvious way to try to implement a controller in our simulation would be to track the change in its position from frame to frame, convert that into an isometry, and move the virtual controller by that isometry. However, this would run into problems with parallel transport, similar to the floor falling away from the user in $\HH^3$: as you walk forward, your hand would appear to diverge from your path, sliding off into the distance. Instead, we plan to update the position of the controller each frame as an offset isometry from the position of the headset. With the correct choice of scaling between real-life euclidean space and our virtual space, this should mean that the controller appears in a location consistent with the user's sense of proprioception.

Floating point errors provide a challenge in tracking objects in hyperbolic space, since the elements of our isometry matrices are exponential functions of the distance travelled in applying them. In our simulations so far, the user never actually leaves the central cube of the tiling: as they attempt to they are teleported from one side of it to the other, and the colours of the cells updated appropriately. We would not be able to use this trick for objects left in the world. One solution would be to record the location of an object by its position relative to a cube of our tiling, together with the cube's location recorded as a word in the tiling's generators.

\end{document}